\documentclass[11pt]{article}
\usepackage{amsfonts}
\usepackage{bbm}
\usepackage{amscd}
\usepackage{xypic}
\usepackage{amsmath}
\usepackage{amssymb}
\usepackage{amsthm}
\usepackage{bbm}
\usepackage{CJK}
\usepackage{fancyhdr}
\usepackage{hyperref}
\usepackage{indentfirst}
\usepackage{latexsym}
\usepackage{mathrsfs}
\allowdisplaybreaks 
\usepackage{color}
\usepackage{amsmath}
\numberwithin{equation}{section}
\usepackage[top=1in,bottom=1in,left=1.25in,right=1.25in]{geometry}
\def\<{\langle}
\def\>{\rangle}

\def\<{\langle}
\def\>{\rangle}

\date{}

\begin{document}
	\renewcommand{\baselinestretch}{1.2}
	\renewcommand{\arraystretch}{1.0}
	\title{\bf Rota-Baxter type $H$-operators on pseudoalgebras }
	\date{}
\author {{\bf Botong Gai}$^a$\,and \, {\bf Shuanhong Wang$^b$\footnote {Corresponding author:  shuanhwang@seu.edu.cn}}\\
{\small a: School of Mathematics, Southeast University,  Nanjing, Jiangsu}\\
{\small 210096, P. R. of China. E-mail: 230228425@seu.edu.cn}\\
{\small b:  Shing-Tung Yau Center, School of Mathematics, Southeast University}\\	
{\small Nanjing, Jiangsu 210096, P. R. of China}}
 \maketitle
	\begin{center}
		\begin{minipage}{14.cm}
			{\bf Abstract} Let $H$ be a Hopf algebra. In this paper, we study a class of $H$-operators on $H$-pseudoalgebras, which
resemble the Rota-Baxter $H$-operator, and they are called Rota-Baxter type $H$-operators. We firstly present some basic properties and examples. Then by using Rota-Baxter type $H$-operators, we construct a number of associative (resp. Lie, NS-) $H$-pseudoalgebras. Meanwhile, Rota-Baxter type $H$-operators on $H$-pseudoalgebras of rank one are studies respectively. Finally, we consider the annihilation algebras and $H$-conformal algebras induced by $H$-pseudoalgebras and corresponding Rota-Baxter type operators are discussed.
\\

\\

		\end{minipage}
	\end{center}

\section*{Introduction}
	\def\theequation{0. \arabic{equation}}
	\setcounter{equation} {0} \hskip\parindent

The notion of conformal algebra (\cite{K}) was introduced by Kac as an axiomatic description of
the operator product expansion (OPE) of chiral fields in conformal field theory, and it came to
be useful for investigation of vertex algebras. Recall that a Lie conformal algebra $L$ is defined
as a $\mathbb{C}[\partial]$-module ($\partial$ is an indeterminate), endowed with a $\mathbb{C}$-linear map
\begin{equation*}
  L\otimes L\rightarrow \mathbb{C}[\lambda]\otimes L, \quad a\otimes b\rightarrow [a_{\lambda}b]
\end{equation*}
satisfying axioms similar to those of Lie algebra (\cite{DK}).
\\

 In 2001, Bakalov, D'Andrea and Kac introduced the notion of Lie $H$-pseudoalgebras by replacing the above polynomial algebra $\mathbb{C}[\partial]$
  with any cocommutative Hopf algebra $H$ (see \cite{BDK}). Furthermore,  they established such theory as
  cohomology theory, representation theory of $H$-pseudoalgebras
   and  irreducible modules over finite simple Lie pseudoalgebras (see \cite{BDK}-\cite{BDK3}, \cite{DL}). Actually,
  Lie $H$-pseudoalgebras  can be considered as Lie algebras in a certain ``pseudotensor'' category, instead of the category
of vector spaces. A pseudotensor category (\cite{BD}) is a category equipped with ``polylinear maps'' and a
way to compose them (such categories were first introduced by Lambek (\cite{L}) under the name multi-categories). This is enough to define the notions of Lie algebra, representations, cohomology, etc.
  Besides, (Lie) $H$-pseudoalgebras are closely related to the differential Lie algebras of the Ritt and Hamiltonian formalism in the theory of nonlinear evolution equations (see \cite{D}\cite{G1}\cite{G2}). However, their roles in many fields of mathematical physics are not yet completely understood since they are relatively new algebraic structures.  Recently there has been many interest in the theory of $H$-pseudoalgebras (for example, \cite{BL}\cite{SQX}\cite{SW}).
\\

The notion of Rota-Baxter operator, introduced by G. Baxter, has been intensively studied in various
fields of mathematics, where mainly its relations with combinatorics and many
other branches of mathematical physics were studied (see \cite{GK},\cite{RK}). For a $\mathbf{k}$-algebra $A$ and a scalar $\lambda\in\mathbf{k}$, a linear operator $R: A\rightarrow A$ is called an Rota-Baxter operator of weight $\lambda$ if the following equation holds for any $a, b\in A$:
 \begin{equation*}
   R(a)R(b)=R[R(a)b+aR(b)+\lambda ab].
 \end{equation*}
 One can refer to the book \cite{G} for the
concrete theory of Rota-Baxter algebras. In recent years, many linear operators have been interested and studied. Examples are:
\begin{align*}
  Endomorphism~operator \qquad &P(ab)=P(a)P(b),\\
  Differential~operator \qquad &P(ab)=P(a)b+aP(b),\\
  Averaging~operator \qquad & P(a)P(b)=P(P(a)b)=P(aP(b)),\\
  Nijenhuis~operator \qquad &P(a)P(b)=P[P(a)b+aP(b)-P(ab)],\\
  Reynolds~operator~of~weight~\lambda \qquad &P(a)P(b)=P[P(a)b+aP(b)+\lambda P(a)P(b)].
\end{align*}
These operators can be divided into two types. The first two operators have the form $P(ab)=N(a, b)$, where $N(a, b)$ is some algebraic expression involving $a, b$ and the operator $P$, and this type is called differential type operators in many articles because of their resemblance to the differential operator. The remaining operators satisfy
an identity of the form $P(a)P(b)=P[N(a, b)]$, where $N(a, b)$ is some algebraic expression involving $a, b$ and the operator $P$, and this type is called Rota-Baxter type operators because of their resemblance to the Rota-Baxter operator.
\\

Averaging operators first appeared in the work of Reynolds in 1895 in connection with the theory of turbulence \cite{R0}.  In a series of papers of 1930s, Kolmogoroff and Kamp'e de F'eriet introduced explicitly the averaging operator in the context of turbulence theory and functional
analysis \cite{KF}\cite{M}. Birkhoff continued the line of research in \cite{B}. Moy investigated averaging operators from the viewpoint of conditional expectation in probability theory \cite{MS}. Later, Kelley \cite{KJL} and Rota \cite{RGC} studied the role of averaging operators in Banach algebras. Contrary to the above studies of analytic nature, the algebraic study on averaging operators began with the Ph.D thesis of Cao \cite{C}. In paper \cite{PG}, Guo and Pei studied averaging operators from an algebraic and combinatorial point
of view.
\\

The study of Nijenhuis operators dated back to late 1970s. It was
discovered by Gel'fand and Dorfman (\cite{GY1},\cite{GY2}) that Nijenhuis operators are closely related
to Hamiltonian pairs. Interrelations between Nijenhuis operators and deformations of Lie
algebras were presented by Dorfman (\cite{D2}). Carinena and his coauthors (\cite{CGM}) introduced an
associative version of classical Nijenhuis identity. For more studies about Nijenhuis operators, see \cite{LP} and \cite{DS}.\\

Reynolds operators were introduced by Reynolds in \cite{R} in the study of fluctuation theory in fluid dynamics. In \cite{F}, the author coined the concept of the Reynolds operator and regarded the operator as a mathematical subject in general.
The authors provided examples and properties of Reynolds operators and studied
the free Reynolds algebras in \cite{GGL} and \cite{ZGG}. More detail of Reynolds operators, see \cite{FM}, \cite{M}, \cite{MM}.\\

The article is organized as follows. In Section 2, we introduce the notions of Rota-Baxter type $H$-operators, which mainly contains averaging $H$-operator, Nijenhuis $H$-operator and Reynolds $H$-operator. Some properties and concrete examples are also given. In section 3, we construct the associative (Lie, NS-) $H$-pseudoalgebras from different Rota-Baxter type $H$-operators. In section 4, we consider the Rota-Baxter type $H$-operators on rank one Lie $H$-pseudoalgebra and present a clear classification. In section 5, we focus on the annihilation algebra of $H$-pseudoalgebra. A relationship between Rota-Baxter type $H$-operators on $H$-pseudoalgebra and Rota-Baxter operators on its annihilation algebra is given. In Section 6, we introduce the conformal averaging (Nijenhuis, Reynolds) $H$-operators on Lie $H$-conformal algebra and discuss the relationship with corresponding $H$-operators on Lie $H$-pseudoalgebra.
\\

Unless otherwise specified, all vector spaces, linear maps and tensor products are considered over an algebraically
 closed field $\mathbf{k}$ of characteristic $0$.\\

\section*{1. Preliminaries}
	\def\theequation{1. \arabic{equation}}
	\setcounter{equation} {0} \hskip\parindent

\subsection*{1.1 Hopf algebra}

In this section we present some facts and notations of Hopf algebra which will be used throughout this paper and can be found, for example, in Sweedler's book (cf. \cite{S}). Let $H$ be a Hopf algebra with a coproduct $\Delta$ , a counit $\epsilon$, and an antipode $S$.
We will use the standard Sweedler's notation (cf. \cite{S}): $\Delta(h) = h_{(1)}\otimes h_{(2)}$. The axioms of the antipode and the counit
can be written as follows:
\begin{equation*}
   S(h_{(1)})h_{(2)}=\epsilon(h)1_{H}=h_{(1)}S(h_{(2)}), \quad \epsilon(h_{(1)})h_{(2)}=h_{(1)}\epsilon(h_{(2)})=h,
\end{equation*}
while the fact that $\Delta$ is an algebra homomorphism translates as
  \begin{equation*}
  (hg)_{(1)}\otimes(hg)_{(2)}=h_{(1)}g_{(1)}\otimes h_{(2)}g_{(2)}.
  \end{equation*}

An element $c\in H$ is called \emph{group-like} if $\Delta(c)=c\otimes c$ and $\epsilon(c)=1$. The set of group-like elements in $H$ is denoted by $G(H)$.

An element $t\in H$ is called \emph{left integral} in $H$ if $ht=\epsilon(h)t$ for all $h\in H$. A \emph{right integral} in $H$ is an element $t^{'}\in H$ such that $t^{'}h=\epsilon(h)t^{'}$ for all $h\in H$.\\

 Since we shall work with cocommutative Hopf algebras, we recall the following important and
 classical result:
 \\

 \textbf{Theorem 1.1 }(\cite{S}, Section 13.1)
 Let $H$ be a cocommutative Hopf algebra over $\mathbf{k}$. Then $H$ is isomorphic (as a Hopf algebra) to the smash product of
 the universal enveloping algebra
$U(P(H))$ and the group algebra $k[G(H)]$, where $G(H)=\{g \in H \mid \Delta(g) = g \otimes g\}$ is the subset of group-like
elements of $H$, and $P(H)=\{p \in H | \Delta(p) = p \otimes 1 + 1 \otimes p\}$ is the subspace of primitive elements of $H$.
\\

For a Hopf algebra $H$, there exists an increasing sequence of subspaces defined inductively by:
\begin{equation*}
  F^{0}H=\mathbf{k}[G(H)],  \quad F^{n}H=0, \quad for ~ n<0,
\end{equation*}
\begin{equation*}
  F^{n}H=span_{\mathbf{k}}\{h\in H|\Delta(h)\in F^{0}H\otimes h+h\otimes F^{0}H+\sum_{i=1}^{n-1}F^{i}H\otimes F^{n-i}H\},~for~n\geq 1.
\end{equation*}
It has the following properties (which are immediate from definitions):
\begin{equation}
  (F^{n}H)(F^{m}H)\subset F^{n+m}H, \quad \Delta(F^{n}H)\subset\sum_{i=1}^{n-1}F^{i}H\otimes F^{n-i}H, \quad S(F^{n}H)\subset F^{n}H.
\end{equation}
When $H$ is cocommutative, one can show that:
\begin{equation}
  \bigcup_{n} F^{n}H=H.
\end{equation}
At this time, we say that a nonzero element $a\in H$ has \textbf{degree $n$} if $a \in F^{n} H \setminus F^{n-1}H$.
\\

When $H$ is a universal enveloping algebra of a finite-dimensional Lie algebra, or its
smash product with the group algebra of a finite group, the following finiteness condition holds:
\begin{equation}
  dim F^{n}H<\infty.
\end{equation}

Further, let $\delta$ be a finite-dimensional Lie algebra and $\{\partial_{i},...,\partial_{N}\}$ be a basis of $\delta$. For $I=(i_{1},i_{2},...,i_{N})\in \mathbb{Z}_{+}^{N}$, let $\partial^{I}=\frac{\partial_{1}^{i_{1}}...\partial_{N}^{i_{N}}}{i_{1}!...i_{N}!}$, then $\{\partial^{I}\}$ is a basis of $H$ (the \emph{PBW-basis}). Moreover, the coproduct on $\partial^{I}$ is
\begin{equation}
  \Delta(\partial^{I})=\sum_{J+K=I}\partial^{J}\otimes\partial^{K}.
\end{equation}
\\

Now let $X=H^{\ast}=Hom_{\mathbf{k}}(H,\mathbf{k})$ be the dual of $H$. The left and right actions of $H$ on $X$ are ($\forall h, f\in H, x\in X$):
\begin{equation*}
  \langle h\cdot x, f \rangle =\langle x, fh \rangle, \qquad \langle x \cdot h, f \rangle =\langle x, hf \rangle.
\end{equation*}
Actually, the associativity of $H$ implies that $X$ is an $H$-bimodule:
\begin{equation*}
  f\cdot(x\cdot g)=(f\cdot x)\cdot g,   \qquad f,g\in H, x\in X.
\end{equation*}
Besides, we can defined an antipode $S^{\ast}: X\rightarrow X$ as the dual of that of $H$:
\begin{equation*}
  \langle S^{\ast}(x), h \rangle=\langle x, S(h)\rangle,
\end{equation*}
then we have $S^{\ast}(xy)=S^{\ast}(y)S^{\ast}(x)$.

We will also define a comultiplication $\Delta: X\rightarrow X\otimes X$ writing as $\Delta(x)=x_{(1)}\otimes x_{(2)}$. By definition, for $x, y\in X, h, g\in H$, we have
\begin{align}
  &\langle xy, h\rangle=\langle x, h_{(1)}\rangle\langle y, h_{(2)}\rangle,\\
  &\langle x, gh\rangle=\langle x_{(1)}, g\rangle\langle x_{(2)}, h\rangle.
\end{align}
\\

 For an arbitrary Hopf algebra $H$,  a map $\mathcal{F} : H \otimes H \rightarrow H \otimes H$ is called the
\emph{Fourier transform} (or \emph{Galois map}), by the formula:
\begin{equation*}
  \mathcal{F}(f\otimes g)=fS(g_{(1)})\otimes g_{(2)},
\end{equation*}
for all $f, g\in H$. Observe that $\mathcal{F}$ is a vector space isomorphism with an inverse given by
\begin{center}
  $\mathcal{F}^{-1}(f\otimes g)=fg_{(1)}\otimes g_{(2)}$.
\end{center}
The significance of $\mathcal{F}$ is in the identity
\begin{equation}
  f\otimes g=\mathcal{F}^{-1}\mathcal{F}(f\otimes g)=(fS(g_{(1)})\otimes 1)\Delta(g_{(2)}),
\end{equation}
which implies the following result.
\\

\textbf{Lemma 1.2.}(\cite{BDK} or \cite{VD1}) Every element of $H \otimes H$ can be uniquely represented in the form $\sum_{i}(h_{i}\otimes 1)\Delta(l_{i})$, where ${h_{i}}$ is a fixed $\mathbf{k}$-basis of $H$ and $l_{i}\in H$. In other words, $H\otimes H = (H \otimes \mathbf{k})\Delta(H)$.
\\

\subsection*{1.2. (Lie) $H$-pseudoalgebras}

In this section, we shall recall the basic definition and example of (Lie) $H$-pseudoalgebra.
\\

\textbf{Definition 1.3}(\cite{BDK}) Let $H$ be an cocommutative Hopf algebra. An (associative) \emph{$H$-pseudoalgebra} $(A, \ast)$ is a left $H$-module $A$ together with a map
(called \emph{pseudoproduct}):
\begin{equation*}
  \ast: A\otimes A\rightarrow (H\otimes H)\otimes_{H} A, \quad a\otimes b\mapsto a\ast b,
\end{equation*}
satisfying the following properties: for any $a, b, c \in L, h, g \in H$,

\textbf{$H$-bilinearity:}
\begin{equation}
  ha\ast gb=(h\otimes g\otimes_{H} 1)(a\ast b),
\end{equation}

\textbf{Associativity:}
\begin{equation}
  a\ast(b\ast c)=(a\ast b)\ast c \in H^{\otimes3}\otimes_{H} A.
\end{equation}

Moreover, for any $\beta\in H\otimes H$, the pseudoproduct is extended by
\begin{align*}
  &(\beta\otimes_{H} a)\ast b=\sum_{i} (\beta\otimes 1)(\Delta\otimes \mathrm{id})(\alpha_{i})\otimes_{H}c_{i}\in H^{\otimes3}\otimes_{H} A,\\
  &a\ast (\beta\otimes_{H} b)=\sum_{i} (1\otimes \beta)(\mathrm{id}\otimes\Delta)(\alpha_{i})\otimes_{H}c_{i}\in H^{\otimes3}\otimes_{H} A,
\end{align*}
if $a\ast b=\sum_{i} \alpha_{i}\otimes_{H} c_{i}\in H\otimes H\otimes_{H} A$.

The pseudoproduct $a\ast b$ is \emph{commutative} iff
\begin{equation*}
  b\ast a=(\sigma\otimes_{H} \mathrm{id})(a\ast b),
\end{equation*}
where $\sigma: H\otimes H$ is the permutation $\sigma(f\otimes g)=g\otimes f$. The right-hand side of above is well
defined due to the cocommutativity of $H$.
\\

Similarly, an \emph{Lie $H$-pseudoalgebra} $(L, [\ast])$ is a left $H$-module $L$ together with a pseudobracket:
\begin{equation*}
  [\ast]: A\otimes A\rightarrow (H\otimes H)\otimes_{H} A, \quad a\otimes b\mapsto [a\ast b],
\end{equation*}
satisfies

\textbf{$H$-bilinearity:}
\begin{equation}
  [ha\ast gb]=(h\otimes g\otimes_{H} 1)[a\ast b],
\end{equation}

\textbf{Skew-commutativity:}
\begin{equation}
  [b\ast a]=-(\sigma\otimes_{H} \mathrm{id})[a\ast b],
\end{equation}

\textbf{Jacobi identity:}
\begin{equation}
  [a\ast[b\ast c]]=[[a\ast b]\ast c]+((\sigma\otimes \mathrm{id})\otimes_{H} \mathrm{id})[b\ast[a\ast c]]\in H^{\otimes3}\otimes_{H} L.
\end{equation}

It is important to point out that the above pseudo-bracket are well defined, provided that the Hopf algebra $H$ is cocommutative.  We will always assume that $H$ is cocommutative when talking about Lie $H$-pseudoalgebras.\\

\textbf{Example 1.4} (\emph{Current pseudoalgebras}) Let $H'$ be a Hopf subalgebra of $H$ and $\mathcal{A}$ be a $H'$-pseudoalgebra. Then we define the current $H$-pseudoalgebra $Cur^{H}_{H'} \mathcal{A} \equiv Cur \mathcal{A}$ as $H\otimes_{H'} \mathcal{A}$ by extending the
pseudoproduct $a\ast b$ of $\mathcal{A}$ using the $H$-bilinearity. Explicitly, for $a, b\in \mathcal{A}$, $f, g\in H$,
\begin{center}
  $(f\otimes_{H'} a) \ast (g\otimes_{H'} b)=(f\otimes g\otimes_{H} 1)(a\ast b)$.
\end{center}

An important special case is when $H'=\mathbf{k}$: given a algebra $A$, let
$Cur A=H\otimes A$ with the following pseudoproduct:
\begin{equation*}
  (f\otimes a) \ast (g\otimes b)=(f\otimes g)\otimes_{H} (1\otimes ab).
\end{equation*}
Then $Cur A$ is a $H$-pseudoalgebra.
\\

Let $A_{1}, A_{2}$ be $H$-pseudoalgebras. A \emph{pseudoalgebra homomorphism} from $A_{1}$ to $A_{2}$ is an $H$-linear map $\varphi: A_{1}\rightarrow A_{2}$ such that
\begin{equation}
  \varphi(a)\ast\varphi(b)=(\mathrm{id}\otimes_{H}\varphi)(a\ast b), \quad \forall a, b\in A.
\end{equation}

Let $A$ be an $H$-pseudoalgebra, $a, b\in A$ and $\{h^{i}\}$ be a $\mathbf{k}$-basis of $H$. By Lemma 1.2, $a\ast b$ can be uniquely written in the form
\begin{equation*}
  a\ast b=\sum_{i}h^{i}\otimes1\otimes_{H}c_{i}.
\end{equation*}
We will call the elements $c_{i}\in A$ \emph{coefficients} of $ a\ast b$. A \emph{subalgebra} of pseudoalgebra $A$ is an $H$-submodule $M$ such that for any $a, b\in M$, all coefficients of $a\ast b$ lie in $M$. An \emph{ideal} of pseudoalgebra $A$ is an $H$-submodule $I$ such that for any $a\in A, b\in I$, all coefficients of $a\ast b$ lie in $I$. A \emph{representation} of an associative $H$-pseudoalgebra $A$, or an \emph{$A$-module}, is a left $H$-module $M$ endowed with a pseudoaction $\ast$ satisfying:
\begin{equation*}
  (a\ast b)\ast m=a\ast(b\ast m)
\end{equation*}
for any $a, b\in A, m\in M$.

Let $M_{1}, M_{2}$ be $A$-modules. An \emph{$A$-module homomorphism} from $M_{1}$ to $M_{2}$ is
an $H$-homomorphism $\rho: M_{1}\rightarrow M_{2}$ such that for any $a\in A, m\in M_{1}$,
\begin{equation}
  a\ast\rho(m)=(\mathrm{id}\otimes_{H} \rho)(a\ast m).
\end{equation}
\\

For a  cocommutative Hopf algebra $H$ and a group $\Gamma$ which acts on $H$ by automorphisms, there is the smash product $\widetilde{H}=H\sharp \mathbf{k}[\Gamma] $of
$H$ with the group algebra of $\Gamma$, and we have the following results.\\

\textbf{Theorem 1.5}(\cite{BDK}) A Lie $\widetilde{H}=H\sharp \mathbf{k}[\Gamma]$-pseudoalgebra $L$ is the same as
a Lie $H$-pseudoalgebra $L$ on which the group $\Gamma$ acts in a way compatible with
the action of $H$, by preserving the $H$-pseudobracket:
\begin{center}
  $[ga\ast gb]=g\cdot [a\ast b]$ \quad for $g\in \Gamma$, $a, b\in L$,
\end{center}
and satisfying the following finiteness condition:
\begin{center}
  given $a, b \in L$, $[ ga \ast b]\neq 0$  for only a finite number of $g \in \Gamma$.
\end{center}
The $\widetilde{H}$-pseudobracket of $L$ is given by the formula
\begin{center}
  $[a \widetilde{\ast} b]=\sum_{g\in\Gamma}((g^{-1}\otimes 1)\otimes_{\widetilde{H}} 1)[ga\ast b]$ \quad for $a, b\in L$.
\end{center}

\textbf{Remark:} This result, combined with Theorem 1.1, will allow us in many
cases to reduce the study of Lie $H$-pseudoalgebras to the case when $H$ is an universal enveloping algebra.\\

\section*{2. Rota-Baxter type $H$-operators}
\def\theequation{2. \arabic{equation}}
	\setcounter{equation} {0} \hskip\parindent

Rota-Baxter $H$-operator on (Lie) $H$-pseudoalgebras was first introduced in \cite{LW}, which is an $H$-linear map $R$ on $H$-pseudoalgebra $A$ and satisfies
\begin{equation*}
  R(a)\ast R(b)=(\mathrm{id}\otimes_{H} R)(R(a)\ast b+a\ast R(b)+\lambda a\ast b)
\end{equation*}
for any $a, b\in A$.

In this section, we introduce a class of $H$-operators on (Lie) $H$-pseudoalgebras, which mainly contains averaging $H$-operator,
Nijenhuis $H$-operator and Reynolds $H$-operator. These operators are called Rota-Baxter type $H$-operators since they resemble the Rota-Baxter $H$-operator.
 \\

\subsection*{2.1 Averaging $H$-operator}

\textbf{Definition 2.1} Let $(A, \ast)$ be a $H$-pseudoalgebra. An $H$-linear map $T: A\rightarrow A$ is called an averaging $H$-operator if
\begin{equation}
  (\mathrm{id}\otimes_{H} T)(T(a)\ast b)=T(a)\ast T(b)=(\mathrm{id}\otimes_{H} T)(a\ast T(b))
\end{equation}
satisfied for all $a, b\in A$.

An $H$-pseudoalgebra $(A, \ast)$ with an averaging $H$-operator $T$
is called an averaging $H$-pseudoalgebra.\\

\textbf{Remark.} An averaging $H$-pseudoalgebra is exactly an averaging algebra
when $H=\mathbf{k}$.\\

\textbf{Example 2.2} (1). Let $A$ be an $H$-pseudoalgebra, then $T:A\rightarrow A$, $a\mapsto ka$ is an averaging $H$-operator for any $k\in
\mathbf{k}$.

(2). Let $A=He_{1}\oplus He_{2}$ be an $H$-pseudoalgebra of rank two with unique nonzero pseudoproduct $e_{2}\ast e_{2}=\alpha\otimes_{H} e_{2}$ for nonzero $\alpha\in H\otimes H$.
Define
\begin{equation*}
  T(e_{1})=he_{1}, \qquad T(e_{2})=\lambda e_{2},
\end{equation*}
then $(A, \ast, T)$ is an averaging $H$-pseudoalgebra for any $h\in H$, $\lambda\in\mathbf{k}$.

(3). Let $A$ be a $\mathbf{k}$-algebra and $T$ an averaging operator on $A$, then $\mathrm{id}_{H}\otimes T$ is an averaging $H$-operator on current
pseudoalgebra $\emph{Cur}A=H\otimes A$.\\

As an $H$-linear map on pseudoalgebras, averaging $H$-operator $T$ has the following properties.\\

\textbf{Proposition 2.3} Let $(A, \ast, T)$ be an averaging $H$-pseudoalgebra, then

(1). $T^{n}(A)$ is a subalgebra of $A$,

(2). $A$ is left $T(A)$-module with left action pseudoproduct,

(3). $T: A\rightarrow A$ is left $T(A)$-module map.

\textbf{Proof.} (1). We only need to show $T(A)$ is a subalgebra of $A$. Firstly, $T(A)$ is an $H$-submodule directly comes from the $H$-linearity of $T$. Besides, from (2.1), $T(A)\ast T(A)$ must lie in $H\otimes H\otimes_{H} T(A)$, which shows $T(A)$ is a subalgebra.

(2). $A$ is left $T(A)$-module directly comes from the associativity of pseudoproduct.

(3). For any $T(a)\in T(A), b\in A$, we have
\begin{equation*}
T(a)\ast T(b)\overset{(2.1)}{=}(\mathrm{id}\otimes T)(T(a)\ast b),
\end{equation*}
which shows $T$ is a left $T(A)$-module map from (1.14).
$\hfill \blacksquare$
\\

\textbf{Proposition 2.4} Let $(A, \ast, T)$ be an averaging $H$-pseudoalgebra, then

(1). $(A, \ast, kT)$ is an averaging $H$-pseudoalgebra for any $k\in \mathbf{k}$,

(2). $T^{\tau}=\tau^{-1} T\tau$ is also an averaging $H$-operator on $A$ where $\tau: A\rightarrow A$ is an automorphism of $A$.

\textbf{Proof.} (1). It can be verified by a straightforward computation.

(2). Since $\tau$ is an automorphism, we have
\begin{equation}
  \tau(a)\ast\tau(b)=(\mathrm{id}\otimes_{H} \tau)(a\ast b)
\end{equation}
for any $a, b\in A$.
Set $u=\tau(a), v=\tau(b)$, then $a=\tau^{-1}(u), b=\tau^{-1}(v)$ and from (2.2) we obtain
\begin{equation}
  \tau^{-1}(u)\ast\tau^{-1}(v)=(\mathrm{id}\otimes_{H} \tau^{-1})(u\ast v).
\end{equation}
Therefore,
\begin{align*}
  (\mathrm{id}\otimes T^{\tau})(T^{\tau}(a)\ast b)=&(\mathrm{id}\otimes_{H} \tau^{-1})(\mathrm{id}\otimes_{H} T)(\mathrm{id}\otimes_{H} \tau)(\tau^{-1} T\tau(a)\ast b) \\
  \overset{(2.2)}{=}&(\mathrm{id}\otimes_{H} \tau^{-1})(\mathrm{id}\otimes_{H} T)(T\tau(a)\ast \tau(b))\\
   \overset{(2.1)}{=}&(\mathrm{id}\otimes_{H} \tau^{-1})(T\tau(a)\ast T\tau(b))\\
   \overset{(2.3)}{=}&T^{\tau}(a)\ast T^{\tau}(b).
\end{align*}
The other side is similar.
$\hfill \blacksquare$
\\

\textbf{Proposition 2.5} Let $T_{1}$ and $T_{2}$ be averaging $H$-operators on $H$-pseudoalgebra $A$, if $T_{1}T_{2}=T_{2}T_{1}$, then $T_{1}T_{2}$ is also an averaging $H$-operator on
 $A$. In particular, a power of an averaging $H$-operator is also an averaging operator.

 \textbf{Proof.} For any $a, b\in A$, we have
 \begin{align*}
     (\mathrm{id}\otimes_{H} T_{1}T_{2})(T_{1}T_{2}(a)\ast b)=&(\mathrm{id}\otimes_{H} T_{2})(\mathrm{id}\otimes_{H} T_{1})(T_{1}[T_{2}(a)]\ast b)\\
     =&(\mathrm{id}\otimes_{H} T_{2})(T_{1}[T_{2}(a)]\ast T_{1}(b))\\
     =&(\mathrm{id}\otimes_{H} T_{2})(T_{2}[T_{1}(a)]\ast T_{1}(b))\\
     =&T_{2}[T_{1}(a)]\ast T_{2}[T_{1}(b)]=T_{1}T_{2}(a)\ast T_{1}T_{2}(b).
 \end{align*}
 The other side is similar.
$\hfill \blacksquare$
\\

\textbf{Proposition 2.6} Let $T_{1}$ and $T_{2}$ be averaging $H$-operators on $H$-pseudoalgebra $A$, then $T_{1}+T_{2}$ is averaging $H$-operators if and only if
\begin{align*}
  T_{1}(a)\ast T_{2}(b)+T_{2}(a)\ast T_{1}(b)=&(\mathrm{id}\otimes_{H} T_{2})(T_{1}(a)\ast b)+(\mathrm{id}\otimes_{H} T_{1})(T_{2}(a)\ast b)\\
  =&(\mathrm{id}\otimes_{H} T_{1})(a\ast T_{2}(b))+(\mathrm{id}\otimes_{H} T_{2})(a\ast T_{1}(b)),
\end{align*}
for any $a, b\in A$.

\textbf{Proof.} On the one hand,
\begin{align*}
  &(\mathrm{id}\otimes_{H}(T_{1}+T_{2}))((T_{1}+T_{2})(a)\ast b)\\
  =&T_{1}(a)\ast T_{1}(b)+T_{2}(a)\ast T_{2}(b)+(\mathrm{id}\otimes_{H} T_{1})(T_{2}(a)\ast b)+(\mathrm{id}\otimes_{H} T_{2})(T_{1}(a)\ast b),
\end{align*}
\begin{align*}
  &(\mathrm{id}\otimes_{H}(T_{1}+T_{2}))(a\ast(T_{1}+T_{2})(b))\\
  =&T_{1}(a)\ast T_{1}(b)+T_{2}(a)\ast T_{2}(b)+(\mathrm{id}\otimes_{H} T_{1})(a\ast T_{2}(b))+(\mathrm{id}\otimes_{H} T_{2})(a\ast T_{1}(b)).
\end{align*}
On the other hand,
\begin{equation*}
  (T_{1}+T_{2})(a)\ast (T_{1}+T_{2})(b)
  =T_{1}(a)\ast T_{1}(b)+T_{2}(a)\ast T_{2}(b)+T_{1}(a)\ast T_{2}(b)+T_{2}(a)\ast T_{1}(b).
\end{equation*}
From (2.1), the condition is sufficient and necessary.
$\hfill \blacksquare$
\\

\textbf{Proposition 2.7} Let $(A, \ast, T)$ be an averaging $H$-pseudoalgebra and $T$ is invertible, then $T^{-1}$ is also an averaging $H$-operator.

\textbf{Proof.} Since $T$ is surjective, for any $v\in A$, there exists $b\in A$ such that $T^{-1}(v)=T(b)$. Then for any $u\in A$,
\begin{equation*}
  u\ast T^{-1}(v)=T(T^{-1}(u))\ast T(b)=(\mathrm{id}\otimes_{H} T)(T^{-1}(u)\ast T(b))=(\mathrm{id}\otimes_{H} T)(T^{-1}(u)\ast T^{-1}(v)).
\end{equation*}
Besides, since $u\ast T^{-1}(v)=(\mathrm{id}\otimes_{H} T)(\mathrm{id}\otimes_{H} T^{-1})(u\ast T^{-1}(v))$ and $T$ is injective, we obtain $T^{-1}(u)\ast T^{-1}(v)=(\mathrm{id}\otimes_{H} T^{-1})(u\ast T^{-1}(v))$. The other side is similar.
$\hfill \blacksquare$
\\

\textbf{Theorem 2.8} Suppose $(A, \ast, T)$ is an averaging $H$-pseudoalgebra and $T$ is inverse, then for
any polynomial $P(t)\in \mathbf{k}[t]$ which does not have a constant term, $P(T)$ is also an
averaging $H$-operator on $A$.

\textbf{Proof.} Let $P_{m}(T)=k_{1}T+...+k_{m}T^{m}$. We know from Proposition 2.4 that $P_{1}(T)=k_{1}T$ is an averaging $H$-operator. Assume it has been established that for $P_{m}(T)$. For any $k_{m+1}\in \mathbf{k}$, we next prove $P_{m+1}(T)=P_{m}(T)+k_{m+1}T^{m+1}$ is also an averaging $H$-operator. Actually, for any $a, b\in A$, we have
\begin{align*}
  (\mathrm{id}\otimes_{H} P_{m}(T))(k_{m+1}T^{m+1}(a)\ast b)=&\sum_{n=1}^{m}(\mathrm{id}\otimes_{H} k_{n}T^{n})(k_{m+1}T^{m+1}(a)\ast b)\\
  =&\sum_{n=1}^{m}k_{n}k_{m+1}(\mathrm{id}\otimes_{H} T^{n})(T^{n}(T^{m+1-n}(a)\ast b))\\
  =&\sum_{n=1}^{m}k_{n}k_{m+1}T^{n}(T^{m+1-n}(a))\ast T^{n}(b)\\
  =&\sum_{n=1}^{m}k_{m+1}T^{m+1}(a)\ast k_{n}T^{n}(b)\\
  =&k_{m+1}T^{m+1}(a)\ast P_{m}(T)(b),
\end{align*}
and
\begin{align*}
 (\mathrm{id}\otimes_{H} k_{m+1}T^{m+1}(T))(P_{m}(T)(a)\ast b)=&\sum_{n=1}^{m}(\mathrm{id}\otimes_{H} k_{m+1}T^{m+1}(T))(k_{n}T^{n}(a)\ast b)\\
 =&\sum_{n=1}^{m}k_{n}k_{m+1}(\mathrm{id}\otimes_{H} T^{m+1})(T^{m+1}(T^{n-m-1}(a)\ast b))\\
 =&\sum_{n=1}^{m}k_{n}k_{m+1}T^{m+1}(T^{n-m-1}(a))\ast T^{m+1}(b)\\
 =&P_{m}(T)(a)\ast k_{m+1}T^{m+1}(T)(b).
\end{align*}
Similarly, one can also obtain
\begin{eqnarray*}
  &(\mathrm{id}\otimes_{H} P_{m}(T))(a\ast k_{m+1}T^{m+1}(b))=P_{m}(T)(a)\ast k_{m+1}T^{m+1}(b),\\
  &(\mathrm{id}\otimes_{H} k_{m+1}T^{m+1}(T))(a\ast P_{m}(T)(b))=k_{m+1}T^{m+1}(T)(a)\ast P_{m}(T)(b).
\end{eqnarray*}
According to Proposition 2.6, $P_{m+1}(T)=P_{m}(T)+k_{m+1}T^{m+1}$ is an averaging $H$-operator.
$\hfill \blacksquare$
\\

\textbf{Remark.} From Proposition 2.6, if $T$ is an averaging $H$-operator, $T+k\mathrm{id}$ is usually not an averaging $H$-operator when $k\neq0$, which is the reason why we give the condition that $P(t)$ does not have a constant term.
\\

\subsection*{2.2 Nijenhuis $H$-operator}

\textbf{Definition 2.9} Let $(A, \ast)$ be an $H$-pseudoalgebra. An $H$-linear map $N: A\rightarrow A$ is called Nijenhuis $H$-operator if
\begin{equation}
  N(a)\ast N(b)=(\mathrm{id}\otimes_{H} N)(N(a)\ast b+a\ast N(b)-(\mathrm{id}\otimes_{H} N)(a\ast b))
\end{equation}
satisfied for any $a, b\in A$.

An $H$-pseudoalgebra $(A, \ast)$ with a Nijenhuis $H$-operator $N$
is called a Nijenhuis $H$-pseudoalgebra.\\

\textbf{Example 2.10} (1). Let $A=He_{1}\oplus He_{2}$ be an $H$-pseudoalgebra of rank two with unique nonzero pseudoproduct $e_{2}\ast e_{2}=\alpha\otimes_{H}e_{2}$ for $\alpha\in H\otimes H$. Define
  \begin{equation*}
    N_{1}(e_{1})=he_{1}, \qquad N_{1}(e_{2})=\lambda e_{2},
  \end{equation*}
or
\begin{equation*}
 N_{2}(e_{1})=\lambda e_{1}, \qquad N_{2}(e_{2})=ge_{1}+\lambda e_{2},
\end{equation*}
for any $h\in H, g\in Z(H)$, $\lambda\in\mathbf{k}$. Then both $N_{1}$ and $N_{2}$ are all Nijenhuis $H$-operators on $A$.

(2). Let $A$ be a $\mathbf{k}$-algebra, $N$ is a Nijenhuis operator on $A$, then $\mathrm{id}_{H}\otimes N$ is a Nijenhuis $H$-operator on current
pseudoalgebra $\emph{Cur}A=H\otimes A$.
\\

\textbf{Proposition 2.11} Let $(A, \ast, N)$ be an Nijenhuis $H$-pseudoalgebra, then

(1). $N^{n}(A)$ is a subalgebra of $A$;

(2). $A$ is left $N(A)$-module with left action pseudoproduct.

\textbf{Proof.} (1). This is similar to Proposition 2.3 by using (2.4).

(2). It is directly comes from the associativity of pseudoproduct.
$\hfill \blacksquare$.
\\

The following Proposition describes close interrelations between Nijenhuis $H$-operators
and Rota-Baxter $H$-operators, and the proof is straightforward.\\

\textbf{Proposition 2.12} Let $N : A \rightarrow A$ be an $H$-module homomorphism over an pseudoalgebra $A$.

(1). If $N^{2}=0$, then $N$ is a Nijenhuis $H$-operator if and only if $N$ is a Rota-Baxter $H$-operator
of weight 0.

(2). If $N^{2}=N$, then $N$ is a Nijenhuis $H$-operator if and only if $N$ is a Rota-Baxter $H$-operator
of weight -1.

(3). If $N^{2}=\mathrm{id}$, then $N$ is a Nijenhuis $H$-operator if and only if $N \pm \mathrm{id}$ is a Rota-Baxter
$H$-operator of weight $\mp2$.
\\

\textbf{Proposition 2.13} Let $(A, \ast, N)$ be a Nijenhuis $H$-pseudoalgebra, then

(1). $(A, \ast, kN)$ is a Nijenhuis $H$-pseudoalgebra for any $k\in \mathbf{k}$;

(2). $N^{\tau}=\tau^{-1} N\tau$ is also a
Nijenhuis $H$-operator if $\tau: A\rightarrow A$ is an automorphism of $A$.

\textbf{Proof.} (1). It can be verified by a straightforward computation.

(2). This is similar to Proposition 2.4.
$\hfill \blacksquare$.
\\

\textbf{Proposition 2.14} Let $(A, \ast, N)$ be an Nijenhuis $H$-pseudoalgebra and $N$ is invertible, then $N^{-1}$ is also an Nijenhuis $H$-operator.

\textbf{Proof.} Since $N$ is surjective, for any $u, v\in A$, there exists $a, b\in A$ such that $N^{-1}(u)=N(a), N^{-1}(v)=N(b)$. Then we have
\begin{align*}
&u\ast N^{-1}(v)+N^{-1}(u)\ast v-(\mathrm{id}\otimes_{H} N^{-1})(u\ast v)\\
=&N^{2}(a)\ast N(b)+N(a)\ast N^{2}(b)-(\mathrm{id}\otimes_{H} N^{-1})[N^{2}(a)\ast N^{2}(b)]\\
\overset{(2.4)}{=}&N^{2}(a)\ast N(b)+N(a)\ast N^{2}(b)-[N^{2}(a)\ast N(b)+N(a)\ast N^{2}(b)-(\mathrm{id}\otimes_{H} N)(N(a)\ast N(b))]\\
 =&(\mathrm{id}\otimes_{H} N)(N(a)\ast N(b))\\
 =&(\mathrm{id}\otimes_{H} N)(N^{-1}(u)\ast N^{-1}(v)).
\end{align*}
Therefore, we obtain $(\mathrm{id}\otimes_{H} N^{-1})[u\ast N^{-1}(v)+N^{-1}(u)\ast v-(\mathrm{id}\otimes N^{-1})(u\ast v)]=N^{-1}(u)\ast N^{-1}(v)$ for any $u, v\in A$.
$\hfill \blacksquare$.
\\

\textbf{Proposition 2.15} Let $N_{1}$ and $N_{2}$ be Nijenhuis $H$-operators on pseudoalgebra $A$, then $N_{1}+N_{2}$ is Nijenhuis $H$-operator if and only if
\begin{align*}
  N_{1}(a)\ast N_{2}(b)+N_{2}(a)\ast N_{1}(b)=&(\mathrm{id}\otimes_{H} N_{2})[N_{1}(a)\ast b+a\ast N_{1}(b)-(\mathrm{id}\otimes_{H} N_{1})(a\ast b)]\\
  +&(\mathrm{id}\otimes_{H} N_{1})[N_{2}(a)\ast b+a\ast N_{2}(b)-(\mathrm{id}\otimes_{H} N_{2})(a\ast b)]
\end{align*}
for any $a, b\in A$.

\textbf{Proof.} This can be easily verified by (2.4).
\\

\textbf{Lemma 2.16} Let $A$ be an $H$-pseudoalgebra and $N: A\rightarrow A$ a Nijenhuis $H$-operator, then for any $a, b\in A$ and any nonnegative numbers $i, j\in \mathbb{Z}$, the following equation holds:
\begin{equation}
  N^{i}(a)\ast N^{j}(b)-(\mathrm{id}\otimes_{H} N^{j})(N^{i}(a)\ast b)-(\mathrm{id}\otimes_{H} N^{i})(a\ast N^{j}(b))+(\mathrm{id}\otimes_{H} N^{i+j})(a\ast b)=0.
\end{equation}

\textbf{Proof.} We prove by induction. It is easy to verify that (2.5) holds for $i=0$ or $j=0$.

For $i=1$, (2.5) is exactly (2.4) when $j=1$. Suppose (2.5) holds for $j$, that is
\begin{equation}
  N(a)\ast N^{j}(b)-(\mathrm{id}\otimes_{H} N^{j})(N(a)\ast b)-(\mathrm{id}\otimes_{H} N)(a\ast N^{j}(b))+(\mathrm{id}\otimes_{H} N^{j+1})(a\ast b)=0.
\end{equation}
Then for $j+1$, we have
\begin{align*}
  &\underline{N(a)\ast N^{j+1}(b)}-(\mathrm{id}\otimes_{H} N^{j+1})(N(a)\ast b)-(\mathrm{id}\otimes_{H} N)(a\ast N^{j+1}(b))+(\mathrm{id}\otimes_{H} N^{j+2})(a\ast b)\\
  \overset{(2.4)}{=}&(\mathrm{id}\otimes_{H} N)[N(a)\ast N^{j}(b)+a\ast N^{j+1}(b)-(\mathrm{id}\otimes_{H} N)(a\ast N^{j}(b))]-(\mathrm{id}\otimes_{H} N^{j+1})(N(a)\ast b)\\
  &-(\mathrm{id}\otimes_{H} N)(a\ast N^{j+1}(b))+(\mathrm{id}\otimes_{H} N^{j+2})(a\ast b)\\
  =&(\mathrm{id}\otimes_{H} N)[N(a)\ast N^{j}(b)+\underline{a\ast N^{j+1}(b)}-(\mathrm{id}\otimes_{H} N)(a\ast N^{j}(b))-(\mathrm{id}\otimes_{H} N^{j})(N(a)\ast b)\\
  &-\underline{a\ast N^{j+1}(b)}+(\mathrm{id}\otimes_{H} N^{j+1})(a\ast b)]\\
  \overset{(2.6)}{=}&0.
\end{align*}
Thus, (2.5) holds for $i=1$ and any $j$.

Next applying (2.6) to the element $N^{i}(a)$ instead of $a$ gives
\begin{equation}
  N^{i+1}(a)\ast N^{j}(b)-(\mathrm{id}\otimes_{H} N^{j})(N^{i+1}(a)\ast b)-(\mathrm{id}\otimes_{H} N)(N^{i}(a)\ast N^{j}(b))+(\mathrm{id}\otimes_{H} N^{j+1})(N^{i}(a)\ast b)=0.
\end{equation}

Suppose (2.5) holds for $i$, then for $i+1$, we have
\begin{align*}
  &\underline{N^{i+1}(a)\ast N^{j}(b)}-(\mathrm{id}\otimes_{H} N^{j})(N^{i+1}(a)\ast b)-(\mathrm{id}\otimes_{H} N^{i+1})(a\ast N^{j}(b))+(\mathrm{id}\otimes_{H} N^{j+i+1})(a\ast b)\\
  \overset{(2.7)}{=}&(\mathrm{id}\otimes_{H} N)(N^{i}(a)\ast N^{j}(b))-(\mathrm{id}\otimes_{H} N^{j+1})(N^{i}(a)\ast b)-(\mathrm{id}\otimes_{H} N^{i+1})(a\ast N^{j}(b))\\
  &+(\mathrm{id}\otimes_{H} N^{j+i+1})(a\ast b)\\
=&(\mathrm{id}\otimes_{H} N)[N^{i}(a)\ast N^{j}(b)-(\mathrm{id}\otimes_{H} N^{j})(N^{i}(a)\ast b)-(\mathrm{id}\otimes_{H} N^{i})(a\ast N^{j}(b))\\
&+(\mathrm{id}\otimes_{H} N^{j+i})(a\ast b)]\\
=&0.
\end{align*}
Thus we have proved the validity of (2.5) for arbitrary $i, j>0$.
$\hfill \blacksquare$.
\\

\textbf{Proposition 2.17} Let $A$ be an $H$-pseudoalgebra and $N : A\rightarrow A$ a Nijenhuis $H$-operator, then for any polynomial $F(t)=\sum_{i=0}^{n}c_{i}t^{i}\in\mathbf{k}[t]$, $F(N)$ is also a Nijenhuis $H$-operator.

\textbf{Proof.} For any $a, b\in A$, we have
\begin{align*}
  &F(N)(a)\ast F(N)(b)-(\mathrm{id}\otimes_{H} F(N))[F(N)(a)\ast b+a\ast F(N)(b)-(\mathrm{id}\otimes_{H}F(N))(a\ast b)]\\
  =&\sum_{i, j=0}^{n}c_{i}c_{j}[N^{i}(a)\ast N^{j}(b)-(\mathrm{id}\otimes_{H} N^{j})(N^{i}(a)\ast b)-(\mathrm{id}\otimes_{H} N^{i})(a\ast N^{j}(b))+(\mathrm{id}\otimes_{H} N^{i+j})(a\ast b)]\\
  \overset{(2.5)}{=}&0.
\end{align*}
The prove is completed.
$\hfill \blacksquare$.
\\

\subsection*{2.3 Reynolds $H$-operator}

\textbf{Definition 2.18} Let $(A, \ast)$ be an $H$-pseudoalgebra. An $H$-linear map $R: A\rightarrow A$ is called a Reynolds $H$-operator of weight $\lambda$ if
\begin{equation}
  R(a)\ast R(b)=(\mathrm{id}\otimes_{H} R)(R(a)\ast b+a\ast R(b)+\lambda R(a)\ast R(b))
\end{equation}
satisfied for any $a, b\in A$.

An $H$-pseudoalgebra $(A, \ast)$ with a Reynolds $H$-operator $R$
is called a Reynolds $H$-pseudoalgebra.\\

\textbf{Remark.} If $R$ is a Reynolds $H$-operator of weight $\lambda$, then $\lambda R$ is a Reynolds $H$-operator of weight 1.
\\

\textbf{Example 2.19} (1). Let $A=He_{1}\oplus He_{2}$ be a pseudoalgebra of rank two with unique nonzero pseudoproduct $e_{2}\ast e_{2}=\alpha\otimes_{H}e_{2}$ for $\alpha\in H\otimes H$. Define                                                                                                                                     \begin{equation*}
  R(e_{1})=he_{1}, \qquad R(e_{2})=\frac{1}{\lambda} e_{2},
\end{equation*}
then $R$ is a Reynolds $H$-operator of weight $\lambda$ on $A$ for any $h\in H$, .

(2). Let $A$ be a $\mathbf{k}$-algebra, $R$ is a Reynolds operator on $A$, then $\mathrm{id}_{H}\otimes R$ is a Reynolds $H$-operator on current
pseudoalgebra $\emph{Cur}A=H\otimes A$.
\\

\textbf{Proposition 2.20} Let $(A, \ast, R)$ be a Reynolds $H$-pseudoalgebra, then

(1). $R^{n}(A)$ is a subalgebra of $A$,

(2). $A$ is left $R(A)$-module with left action pseudoproduct,

\textbf{Proof.} (1). This is similar to Proposition 2.3 by using (2.8).

(2). It is directly comes from the associativity of pseudoproduct.
\\

\textbf{Proposition 2.21} If $(A, \ast, R)$ is an averaging $H$-pseudoalgebra and $R^{2}=R$, then $(A, \ast, R)$ is a Reynolds $H$-pseudoalgebra of weight -1.

\textbf{Proof.} For any $a, b\in A$, we have
\begin{align*}
  &R(a)\ast R(b)+(\mathrm{id}\otimes_{H} R)[R(a)\ast R(b)-a\ast R(b)-R(a)\ast b]\\
  \overset{(2.1)}{=}&R(a)\ast R(b)+(\mathrm{id}\otimes_{H} R)[(\mathrm{id}\otimes_{H} R)(a\ast R(b))-a\ast R(b)]-R(a)\ast R(b)\\
  =&(\mathrm{id}\otimes_{H} R^{2})(a\ast R(b))-(\mathrm{id}\otimes_{H} R)(a\ast R(b))\\
  =&0.\\
\end{align*}
Therefore, $R$ is a Reynolds $H$-operator of weight -1.
$\hfill \blacksquare$.
\\

\textbf{Proposition 2.22} Let $(A, \ast, R)$ be a Reynolds $H$-pseudoalgebra, then

(1). $(A, \ast, kR)$ is a Reynolds $H$-pseudoalgebra for any $k\in \mathbf{k}$,

(2). $R^{\tau}=\tau^{-1} R\tau$ is also a
Reynolds $H$-operator if $\tau: A\rightarrow A$ is an automorphism of $A$.

\textbf{Proof.} (1). It can be verified by a straightforward computation.

(2). This is similar to Proposition 2.4.
$\hfill \blacksquare$.
\\

\textbf{Proposition 2.23} Let $R_{1}$ and $R_{2}$ be Reynolds $H$-operators of weight $\lambda$ on $A$, then $R_{1}+R_{2}$ is Reynolds $H$-operator if and only if
\begin{align*}
  &R_{1}(a)\ast R_{2}(b)+R_{2}(a)\ast R_{1}(b)\\
  =&(\mathrm{id}\otimes_{H} R_{2})[R_{1}(a)\ast b+a\ast R_{1}(b)+\lambda(R_{1}(a)\ast R_{1}(b)+R_{1}(a)\ast R_{2}(b)+R_{2}(a)\ast R_{1}(b))]\\
  +&(\mathrm{id}\otimes_{H} R_{1})[R_{2}(a)\ast b+a\ast R_{2}(b)+\lambda(R_{2}(a)\ast R_{2}(b)+R_{1}(a)\ast R_{2}(b)+R_{2}(a)\ast R_{1}(b))]
\end{align*}
for any $a, b\in A$.

\textbf{Proof.} This can be easily verified by (2.8).
\\

\section*{3. Constructing pseudoalgebras by Rota-Baxter type $H$-operators}
\def\theequation{3. \arabic{equation}}
	\setcounter{equation} {0} \hskip\parindent

In this section, we mainly construct associative (Lie, NS) pseudoalgebras by using Rota-Baxter type $H$-operators over $H$-pseudoalgebras.\\

It has been shown that for an $H$-pseudoalgebra $A$, one can define a pseudobracket
\begin{equation*}
  [a\ast b]=a\ast b-(\sigma\otimes_{H} \mathrm{id})(b\ast a),
\end{equation*}
which makes $A$ a Lie $H$-pseudoalgebra \cite{BDK}. This operation can be
rewritten as
\begin{equation*}
  [a\ast b]=a\ast \mathrm{id}_{A}(b)-(\sigma\otimes_{H} \mathrm{id})(b\ast \mathrm{id}_{A}(a)).
\end{equation*}
This approach can be generalized and make sense by using an averaging $H$-operator on $A$ instead of $\mathrm{id}_{A}$:
\\

\textbf{Proposition 3.1} Let $(A, \ast, T)$ be an averaging $H$-pseudoalgebra, then the $H$-module $A$ will become a Lie $H$-pseudoalgebra under the bracket operation
\begin{equation*}
  [a\ast b]_{T}=a\ast T(b)-(\sigma\otimes_{H} \mathrm{id})(b\ast T(a)).
\end{equation*}

\textbf{Proof.} The $H$-bilinearity of $[\ast]_{T}$ comes from the linearity of $T$. The skew-commutativity is obvious and we only need to verify the Jacobi identity. For any $a, b, c\in A$,
\begin{align*}
  &[[a\ast b]_{T}\ast c]_{T}+((12)\otimes \mathrm{id})[b\ast[a\ast c]_{T}]_{T}-[a\ast[b\ast c]_{T}]_{T}\\
  =&[(a\ast T(b)-(\sigma\otimes_{H} \mathrm{id})(b\ast T(a)))\ast c]_{T}+((12)\otimes_{H} \mathrm{id})[b\ast(a\ast T(c)-(\sigma\otimes_{H} \mathrm{id})(c\ast T(a)))]_{T}\\
  &-[a\ast(b\ast T(c)-(\sigma\otimes_{H} \mathrm{id})(c\ast T(b)))]_{T}\\
  =&a\ast T(b)\ast T(c)-((132)\otimes_{H} \mathrm{id})(c\ast\underline{(\mathrm{id}\otimes_{H} T)(a\ast T(b))})-\{(\sigma\otimes_{H} \mathrm{id})(b\ast T(a))\}\ast T(c)\\
  &+((132)\otimes_{H} \mathrm{id})\{c\ast\underline{(\mathrm{id}\otimes_{H} T)(\sigma\otimes_{H} \mathrm{id})(b\ast T(a))}\}+((12)\otimes_{H} \mathrm{id})\{b\ast\underline{(\mathrm{id}\otimes_{H} T)(a\ast T(c))}\}\\
  &-((12)\otimes_{H} \mathrm{id})((123)\otimes_{H} \mathrm{id})(a\ast T(c)\ast T(b))-((12)\otimes_{H} \mathrm{id})\{b\ast\underline{(\mathrm{id}\otimes_{H} T)(\sigma\otimes_{H} \mathrm{id})(c\ast T(a))}\}\\
  &+((12)\otimes_{H} \mathrm{id})((123)\otimes_{H} \mathrm{id})\{((\sigma\otimes_{H} \mathrm{id})(c\ast T(a)))\ast T(b)\}-a\ast(\underline{(\mathrm{id}\otimes_{H} T)(b\ast T(c))})\\
  &+((123)\otimes_{H} \mathrm{id})(b\ast T(c)\ast T(a))+a\ast\{\underline{(\mathrm{id}\otimes_{H} T)(\sigma\otimes_{H} \mathrm{id})(c\ast T(b))}\}\\
  &-((123)\otimes_{H} \mathrm{id})\{((\sigma\otimes_{H} \mathrm{id})(c\ast T(b)))\ast T(a)\}\\
  \overset{(2.1)}{=}&a\ast T(b)\ast T(c)-((132)\otimes_{H} \mathrm{id})(c\ast T(a)\ast T(b))-\{(\sigma\otimes_{H} \mathrm{id})(b\ast T(a))\}\ast T(c)\\
  &+((132)\otimes_{H} \mathrm{id})\{c\ast((\sigma\otimes_{H} \mathrm{id})(T(b)\ast T(a)))+((12)\otimes_{H} \mathrm{id})\{b\ast T(a)\ast T(c)\}\\
  &-((12)\otimes_{H} \mathrm{id})((123)\otimes_{H} \mathrm{id})(a\ast T(c)\ast T(b))-((12)\otimes_{H} \mathrm{id})\{b\ast((\sigma\otimes_{H} \mathrm{id})(T(c)\ast T(a)))\}\\
  &+((12)\otimes_{H} \mathrm{id})((123)\otimes_{H} \mathrm{id})\{((\sigma\otimes_{H} \mathrm{id})(c\ast T(a)))\ast T(b)\}-a\ast T(b)\ast T(c)\\
  &+((123)\otimes_{H} \mathrm{id})(b\ast T(c)\ast T(a))+a\ast\{(\sigma\otimes_{H} \mathrm{id})(T(c)\ast T(b))\}\\
  &-((123)\otimes_{H} \mathrm{id})\{((\sigma\otimes_{H} \mathrm{id})(c\ast T(b)))\ast T(a)\}\\
  =&0,
\end{align*}
where one can verify terms 1 and 9, 2 and 8, 3 and 5, 4 and 12, 6 and 11, 7 and 10 will cancel each other in  the final step.
$\hfill \blacksquare$
\\

\textbf{Proposition 3.2} Let $(A, \ast, T)$ be an averaging $H$-pseudoalgebra. Define $a\widetilde{\ast}b=a\ast T(b)$ (or $a\widetilde{\ast}b=T(a)\ast b$), then $(A, \widetilde{\ast})$ is an $H$-pseudoalgebra and $T$ is also an averaging $H$-operator on $(A, \widetilde{\ast})$.

\textbf{Proof.} By a straightforward computation,
\begin{center}
$(a\widetilde{\ast}b)\widetilde{\ast}c=(a\ast T(b))\ast T(c)$
\end{center}
and
\begin{center}
$a\widetilde{\ast}(b\widetilde{\ast}c)=a\ast(\mathrm{id}\otimes_{H} T)(b\widetilde{\ast}c)=a\ast(\mathrm{id}\otimes_{H} T)(b\ast T(c))=a\ast(T(b)\ast T(c))$.
\end{center}
Then we omit the remaining details.
$\hfill \blacksquare$.
\\

Nijenhuis operator is closely related to NS-algebra. In the following, we will show a similar result. We firstly give the notion of NS-pseudoalgebra.\\

\textbf{Definition 3.3} Let $\mathcal{A}$ be an $H$-module equipped with three binary operations $\triangleright, \triangleleft, \diamond: \mathcal{A}\otimes \mathcal{A}\rightarrow H^{\otimes2}\otimes_{H} \mathcal{A}$. Then $\mathcal{A}$ is called a \textbf{NS-pseudoalgebra} if $\triangleright, \triangleleft, \diamond$ are $H$-bilinear maps and satisfy the following axioms for any $a, b, c\in \mathcal{A}$:
\begin{align}
&a\triangleright(b\triangleleft c)=(a\triangleright b)\triangleleft c,\\
 &a\triangleright(b\triangleright c)=(a\triangleright b+a \triangleleft b+a\diamond b)\triangleright c,\\
  &a\triangleleft(b\triangleright c+b \triangleleft c+b\diamond c)=(a\triangleleft b)\triangleleft c,\\
a\triangleright(b\diamond c)-(a\triangleright b+&a \triangleleft b+a\diamond b)\diamond c=(a\diamond b)\triangleleft c-a\diamond (b\triangleright c+b \triangleleft c+b\diamond c).
\end{align}
\\

\textbf{Example 3.4} Let $(A, \prec, \succ, \vee)$ be a $NS$-algebra, then $(\emph{CurA}=H\otimes A, \prec_{\ast}, \succ_{\ast}, \vee_{\ast})$ will be a $NS$-pseudoalgebra with $\prec_{\ast}, \succ_{\ast}, \vee_{\ast}$ defined by
\begin{align*}
  (f\otimes a)\prec_{\ast}(g\otimes b)=&(f\otimes g\otimes_{H}(a\prec b))\\
  (f\otimes a)\succ_{\ast}(g\otimes b)=&(f\otimes g\otimes_{H}(a\succ b))\\
  (f\otimes a)\vee_{\ast}(g\otimes b)=&(f\otimes g\otimes_{H}(a\vee b))
\end{align*}
for any $f\otimes a, g\otimes b\in \emph{CurA}$.
\\

\textbf{Proposition 3.5} Let $(\mathcal{N}, \triangleright, \triangleleft, \diamond)$ be a \textbf{NS-pseudoalgebra}, then $(\mathcal{N}, \star)$ forms an associative $H$-pseudoalgebra, where $\star$ is defined by
\begin{equation*}
  a\star b=:a\triangleright b+a \triangleleft b+a\diamond b.
\end{equation*}

\textbf{Proof.} Straightforward.
$\hfill \blacksquare$
\\

The following theorem shows that Nijenhuis $H$-operator on associative $H$-pseudoalgebra
give rise to NS-pseudoalgebra structure.\\

\textbf{Proposition 3.6} Let $(A, \ast, N)$ be a Nijenhuis $H$-pseudoalgebra.
For any $a, b \in A$, define three pseudoproducts on $A$ by
\begin{equation}\label{2}
  a\triangleright b=N(a)\ast b,\quad a\triangleleft b=a\ast N(b), \quad a\diamond b=-(\mathrm{id}\otimes_{H} N)(a\ast b).
\end{equation}
Then $(A, \triangleright, \triangleleft, \diamond)$ is a NS-pseudoalgebra.

\textbf{Proof.} It is easy to see that $\triangleright, \triangleleft$ and $\diamond$ are $H$-bilinear maps. Relation (3.1) follows from the associativity of $A$. Relation (3.2) and (3.3) are actually comes from (2.4). It is left to prove (3.4):
\begin{align*}
  &a\triangleright(b\diamond c)-(a\triangleright b+a \triangleleft b+a\diamond b)\diamond c-(a\diamond b)\triangleleft c+a\diamond (b\triangleright c+b \triangleleft c+b\diamond c)\\
  =&(\mathrm{id}\otimes_{H} N)((N(a)\ast b+a\ast N(b)-(\mathrm{id}\otimes_{H} N)(a\ast b))\ast c)-\underline{N(a)\ast(\mathrm{id}\otimes_{H} N)(b\ast c)}\\
  &+\underline{(\mathrm{id}\otimes_{H} N)(a\ast b)\ast N(c)}-(\mathrm{id}\otimes_{H} N)(a\ast(N(b)\ast c+b\ast N(c)-(\mathrm{id}\otimes_{H} N)(b\ast c)))\\
  \overset{(2.4)}{=}&(\mathrm{id}\otimes_{H} N)[(N(a)\ast b\ast c+a\ast N(b)\ast c-(\mathrm{id}\otimes_{H} N)(a\ast b)\ast c)]-(\mathrm{id}\otimes_{H} N)[a\ast(\mathrm{id}\otimes_{H} N)(b\ast c)\\
  &+N(a)\ast(b\ast c)-(\mathrm{id}\otimes_{H} N)(a\ast b\ast c)]+(\mathrm{id}\otimes_{H} N)[(\mathrm{id}\otimes_{H} N)(a\ast b)\ast c+(a\ast b)\ast N(c)\\
  &-(\mathrm{id}\otimes_{H} N)(a\ast b\ast c)]-(\mathrm{id}\otimes_{H} N)[a\ast N(b)\ast c+a\ast b\ast N(c)-a\ast(\mathrm{id}\otimes_{H} N)(b\ast c)]\\
  =&0.
\end{align*}
This completes the proof.
$\hfill \blacksquare$
\\

\textbf{Proposition 3.7} Let $N$ be a Nijenhuis $H$-operator over a Lie $H$-pseudoalgebra $L$.
Define
\begin{equation}\label{3}
  [a\ast b]_{N}=[N(a)\ast b]+[a\ast N(b)]-(\mathrm{id}\otimes_{H} N)[a\ast b], \quad \forall a, b\in L.
\end{equation}
Then $(L, [\ast]_{N})$ also forms a Lie $H$-pseudoalgebra, denoted by $L^{N}$. Further, $N$ is also
a Lie pseudoalgebra homomorphism from $L^{N}$ to the original Lie pseudoalgebra $L$:
\begin{equation}\label{4}
  [N(a)\ast N(b)]=(\mathrm{id}\otimes_{H} N)[a\ast b]_{N}.
\end{equation}

\textbf{Proof.}
The $H$-bilinearity and Skew-commutativity of $[\ast]_{N}$ comes from the corresponding properties of $[\ast]$. We mainly show the Jacobi identity:
\begin{align*}
  &[[a\ast b]_{N}\ast c]_{N}+((\sigma\otimes\mathrm{id})\otimes_{H} \mathrm{id})[b\ast[a\ast c]_{N}]_{N}-[a\ast[b\ast c]_{N}]_{N}\\
  =&[[N(a)\ast b]\ast c]_{N}+[[a\ast N(b)]\ast c]_{N}-[(\mathrm{id}\otimes_{H} N)[a\ast b]\ast c]_{N}+((\sigma\otimes\mathrm{id})\otimes_{H} \mathrm{id})\{[b\ast[N(a)\ast c]]_{N}\\
&+[b\ast[a\ast N(c)]]_{N}-[b\ast(\mathrm{id}\otimes_{H} N)[a\ast c]]_{N}\}-[a\ast[N(b)\ast c]]_{N}-[a\ast[b\ast N(c)]]_{N}\\&+[a\ast(\mathrm{id}\otimes_{H} N)[b\ast c]]_{N}\\
=&[(\mathrm{id}\otimes_{H} N)[N(a)\ast b]\ast c]+[[N(a)\ast b]\ast N(c)]-(\mathrm{id}\otimes_{H} N)[[N(a)\ast b]\ast c]\\
&+[(\mathrm{id}\otimes_{H} N)[a\ast N(b)]\ast c]+[[a\ast N(b)]\ast N(c)]-(\mathrm{id}\otimes_{H} N)[[a\ast N(b)]\ast c]\\
&-[(\mathrm{id}\otimes_{H} N)^{2}[a\ast b]\ast c]-[(\mathrm{id}\otimes_{H} N)[a\ast b]\ast N(c)]+(\mathrm{id}\otimes_{H} N)[(\mathrm{id}\otimes_{H} N)[a\ast b]\ast c]\\
&+((\sigma\otimes\mathrm{id})\otimes_{H} \mathrm{id})\{[N(b)\ast[N(a)\ast c]]+[b\ast(\mathrm{id}\otimes_{H} N)[N(a)\ast c]]-(\mathrm{id}\otimes_{H} N)[b\ast[N(a)\ast c]]\\
&+[N(b)\ast[a\ast N(c)]]+[b\ast(\mathrm{id}\otimes_{H} N)[a\ast N(c)]]-(\mathrm{id}\otimes_{H} N)[b\ast[a\ast N(c)]]\\&-[N(b)\ast(\mathrm{id}\otimes_{H} N)[a\ast c]]
-[b\ast(\mathrm{id}\otimes_{H} N)^{2}[a\ast c]]+(\mathrm{id}\otimes_{H} N)[b\ast(\mathrm{id}\otimes_{H} N)[a\ast c]]\}\\
&-[N(a)\ast[N(b)\ast c]]-[a\ast(\mathrm{id}\otimes_{H} N)[N(b)\ast c]]+(\mathrm{id}\otimes_{H} N)[a\ast[N(b)\ast c]]-[N(a)\ast[b\ast N(c)]]\\
&-[a\ast(\mathrm{id}\otimes_{H} N)[b\ast N(c)]]+(\mathrm{id}\otimes_{H} N)[a\ast[b\ast N(c)]]+[N(a)\ast(\mathrm{id}\otimes_{H} N)[b\ast c]]\\
&+[a\ast(\mathrm{id}\otimes_{H} N)^{2}[b\ast c]]-(\mathrm{id}\otimes_{H} N)[a\ast(\mathrm{id}\otimes_{H} N)[b\ast c]]\\
=&[[N(a)\ast N(b)]\ast c]+[[N(a)\ast b]\ast N(c)]+[[a\ast N(b)]\ast N(c)]-(\mathrm{id}\otimes_{H} N)[[a\ast b]\ast N(c)]\\
&-(\mathrm{id}\otimes_{H} N)[(\mathrm{id}\otimes_{H} N)[a\ast b]\ast c]+(\mathrm{id}\otimes_{H} N)^{2}[[a\ast b]\ast c]-(\mathrm{id}\otimes_{H} N)[[N(a)\ast b]\ast c]\\
&-(\mathrm{id}\otimes_{H} N)[[a\ast N(b)]\ast c]+(\mathrm{id}\otimes_{H} N)[(\mathrm{id}\otimes_{H} N)[a\ast b]\ast c]\\
&+((\sigma\otimes\mathrm{id})\otimes_{H} \mathrm{id})\{[N(b)\ast[N(a)\ast c]]+[N(b)\ast[a\ast N(c)]]-(\mathrm{id}\otimes_{H} N)[N(b)\ast[a\ast c]]\\
&-(\mathrm{id}\otimes_{H} N)^{2}[b\ast[a\ast c]]+[b\ast[N(a)\ast N(c)]]-(\mathrm{id}\otimes_{H} N)[b\ast[N(a)\ast c]]\\&-(\mathrm{id}\otimes_{H} N)[b\ast[a\ast N(c)]]\}-[N(a)\ast[N(b)\ast c]]-[N(a)\ast[b\ast N(c)]]\\&+(\mathrm{id}\otimes_{H} N)[N(a)\ast[b\ast c]]+(\mathrm{id}\otimes_{H} N)[a\ast(\mathrm{id}\otimes_{H} N)[b\ast c]]-(\mathrm{id}\otimes_{H} N)^{2}[a\ast[b\ast c]]\\&-[a\ast[N(b)\ast N(c)]]+(\mathrm{id}\otimes_{H} N)[a\ast[N(b)\ast c]]+(\mathrm{id}\otimes_{H} N)[a\ast[b\ast N(c)]]\\&-(\mathrm{id}\otimes_{H} N)[a\ast(\mathrm{id}\otimes_{H} N)[b\ast c]]\\
=&0.
\end{align*}
Terms in the penultimate equation can cancel each other by Jacobi identity.
Then $(L, [\ast]_{N})$ is a Lie $H$-pseudoalgebra. Since $N$ is Nijenhuis $H$-operator, then (3.7) naturally holds by (2.4).
$\hfill \blacksquare$
\\

The Rota-Baxter $H$-operator of weight zero is a special case of weight zero Reynolds operator. As in the case of Rota-Baxter pseudoalgebras \cite{LW}, a Reynolds pseudoalgebra structure can replicate itself as follows.\\

\textbf{Proposition 3.8} Let $(A, \ast, R)$ be a Reynolds $H$-pseudoalgebra of weight $\lambda$. Define a new pseudoproduct $\star$ on $A$ by
\begin{equation*}
  a\star b=a\ast R(b)+R(a)\ast b+\lambda R(a)\ast R(b) \qquad \forall a, b\in A,
\end{equation*}
then

(1). $(A, \star)$ is a $H$-pseudoalgebra.

(2). $R$ is also a Reynolds $H$-operator of weight $\lambda$ on $(A, \star)$.

\textbf{Proof.} (1). By a straightforward computation,
\begin{align*}
  &(a\star b)\star c
  =(a\ast R(b)+R(a)\ast b+\lambda R(a)\ast R(b))\star c\\
  =&(a\ast R(b))\ast R(c)+[(\mathrm{id}\otimes_{H} R)(a\ast R(b))]\ast c+[\lambda(\mathrm{id}\otimes_{H} R)(a\ast R(b))]\ast R(c)\\
  &+(R(a)\ast b)\ast R(c)
  +[(\mathrm{id}\otimes_{H} R)(R(a)\ast b)]\ast c+[\lambda(\mathrm{id}\otimes_{H} R)(R(a)\ast b)]\ast R(c)\\
  &+\lambda(R(a)\ast R(b))\ast R(c)
  +[\lambda(\mathrm{id}\otimes_{H} R)(R(a)\ast R(b))]\ast c+\lambda^{2}[(\mathrm{id}\otimes_{H} R)(R(a)\ast R(b))]\ast R(c)\\
 \overset{(2.8)}{=}&R(a)\ast R(b)\ast c+2\lambda R(a)\ast R(b)\ast R(c)+[a\ast R(b)+R(a)\ast b]\ast R(c),
\end{align*}
and
\begin{align*}
  &a\star (b\star c)
  =a\star(b\ast R(c)+R(b)\ast c+\lambda R(b)\ast R(c))\\
  =&a\ast[(\mathrm{id}\otimes_{H} R)(b\ast R(c))]+R(a)\ast b\ast R(c)+\lambda R(a)\ast[(\mathrm{id}\otimes_{H} R)(b\ast R(c))]\\&+a\ast[(\mathrm{id}\otimes_{H} R)(R(b)\ast c)]
  +R(a)\ast R(b)\ast c+\lambda R(a)\ast[(\mathrm{id}\otimes_{H} R)(R(b)\ast c)]\\&+\lambda a\ast[(\mathrm{id}\otimes_{H} R)(R(b)\ast R(c))]
  +\lambda R(a)\ast R(b)\ast R(c)+\lambda^{2} R(a)\ast[(\mathrm{id}\otimes_{H} R)(R(b)\ast R(c))]\\
  \overset{(2.8)}{=}&a\ast R(b)\ast R(c)+2\lambda R(a)\ast R(b)\ast R(c)+R(a)\ast[b\ast R(c)+R(b)\ast c].
\end{align*}
Compare this two results gives ``$\star$'' is associative.

(2). For any $a, b\in A$, we have
\begin{align*}
      &(\mathrm{id}\otimes_{H} R)[a\star R(b)+R(a)\star b+\lambda R(a)\star R(b)]\\
      =&(\mathrm{id}\otimes_{H} R)[R(a)\ast R(b)+R^{2}(a)\ast b+\lambda R^{2}(a)\ast R(b)+
      a\ast R^{2}(b)+R(a)\ast R(b)\\&+\lambda R(a)\ast R^{2}(b)+
      \lambda(R(a)\ast R^{2}(b)+R(a)\ast R(b)+\lambda R^{2}(a)\ast R^{2}(b))]\\
=&R(a)\ast R^{2}(b)+R^{2}(a)\ast R(b)+R^{2}(a)\ast R^{2}(b)\\
=&R(a)\star R(b).
     \end{align*}
     Therefore, $R$ is a Reynolds $H$-operator of weight $\lambda$ on $(A, \star)$.
     $\hfill \blacksquare$
\\

\section*{4. Rota-Baxter type $H$-operators on rank 1 Lie pseudoalgebras}
\def\theequation{4. \arabic{equation}}
	\setcounter{equation} {0} \hskip\parindent

In this section, we mainly discuss Rota-Baxter type $H$-operators on rank 1 Lie pseudoalgebras. Unless otherwise specified, in this section we will always be working with $H=U(\delta)$ where $\delta$ is a finite-dimension Lie algebra.
\\

Let $L=He$ be a Lie $H$-pseudoalgebra. Then, by $H$-bilinearity, the pseudo-bracket on $L$ is determined by $[e\ast e]$, or equivalently, by an $\alpha\in H\otimes H$ such that $[e\ast e]=\alpha\otimes_{H} e$. $\alpha$ is determined by the following Lemma.\\

\textbf{Lemma 4.1}(\cite{BDK}) Let $H=U(\delta)$ be an universal enveloping algebra of a Lie
algebra $\delta$, then $He$ with $[e\ast e]=\alpha\otimes_{H} e$ will be a rank one Lie $H$-pseudoalgebra if and only if $\alpha=r+s\otimes1-1\otimes s$ for some $r\in \delta\wedge\delta, s\in\delta$ satisfy the following system of equations:
\begin{align}
  &[r, \Delta(s)]=0,\\
  &[r_{12}, r_{13}]+r_{12}s_{3}+[r_{12}, r_{23}]+r_{23}s_{1}+[r_{13}, r_{23}]-r_{13}s_{2}=0,
\end{align}
where $r_{12}=r\otimes1, s_{3}=1\otimes1\otimes s$, etc.\\

\textbf{Remark.} For a given Lie algebra $\delta$, a method of searching such $r, s$ in $\delta$ was introduced in \cite{BDK} by using $n$-forms. Actually, if $r=\sum_{i}a_{i}\otimes b_{i}-b_{i}\otimes a_{i}$ and $s$ satisfy (4.1)-(4.2), then there are only two cases need to be considered: when $\delta=\delta_{1}$ and when $\delta=\delta_{1}\oplus\mathbf{k}s$ where $\delta_{1}$ is spanned by (linearly independent set) $\{a_{i}, b_{i}\}$, which give rise to two classes Lie pseudoalgebras: $H(\delta, \chi, \omega)$ and $K(\delta, \theta)$, where $\omega$ (a 2-form), $\chi$ and $\theta$ (1-forms) are determined by $r$ and $s$.

This Lemma actually shows $\alpha$ must lie in $F^{1}(H)\otimes F^{1}(H)$. Based on this result, we can classify all averaging $H$-operator (respectively, Nijenhuis $H$-operator, Reynolds $H$-operator) on $He$.\\

\textbf{Proposition 4.2} Let $H=U(\delta)$ and $L=He$ be a Lie pseudoalgebra of rank 1 with $[e\ast e]=\alpha\otimes_{H} e$ for some nonzero $\alpha$, then the only averaging $H$-operator on $L$ is $T(e)=ke$ for $k\in \mathbf{k}$.

\textbf{Proof.} Suppose $T(e)=he$ for some $h\in H$, then by (2.1), one can obtain
\begin{equation*}
  (h\otimes h)\alpha=(1\otimes h)\alpha\Delta(h)=(h\otimes 1)\alpha\Delta(h).
\end{equation*}
Obviously we have $(1\otimes h-h\otimes1)\alpha\Delta(h)=0$. Since $H$ has no nonzero zerodivisor, then there must have $h\in \mathbf{k}$
and further $h$ is arbitrary.
$\hfill \blacksquare$.
\\

\textbf{Corollary 4.3} Let $H$ be any Hopf algebra with no nonzero divisor and $A=He$ a nontrivial pseudoalgebra of rank 1, then the only averaging $H$-operator on $A$ is $T(e)=ke$ for $k\in\mathbf{k}$.
\\

\textbf{Proposition 4.4} Let $H=U(\delta)$ and $L=He$ be a Lie pseudoalgebra of rank 1 with $[e\ast e]=\alpha\otimes_{H} e$ for some nonzero $\alpha$, then the only Nijenhuis $H$-operator on $L$ is $N(e)=ke$ for $k\in \mathbf{k}$.

\textbf{Proof.} Suppose $N(e)=ge$ for some $g\in H$.
By (2.4), one can obtain
\begin{equation}
  (g\otimes g)\alpha=[(g\otimes1+1\otimes g)\alpha-\alpha\Delta(g)]\Delta(g).
\end{equation}

By Lemma 4.1, $\alpha=r+s\otimes1-1\otimes s\in H\otimes H$.
If $dim\delta>1$, then $r\neq0$.
Set $deg(g)=d$. If $d>1$, then there exist terms in the right hand lie in $F^{2d}(H)\otimes F^{2}(H)$ while all terms in the left hand lie in $F^{d+1}(H)\otimes F^{d+1}(H)$, which is a contradiction. Thus, $d\leq1$ and we can set $g=a+k$ for some $a\in \delta, k\in\mathbf{k}$. Taking this into (4.3), one can obtain
\begin{equation*}
  (a\otimes a)\alpha=[\Delta(a)\alpha, \alpha\Delta(a)]\Delta(a).
\end{equation*}
Specifically,
\begin{equation*}
  (a\otimes a)r+as\otimes a-a\otimes as=[\Delta(a), r]+[a, s]\otimes1-1\otimes[a, s].
\end{equation*}
Notice that there must have $[a, s]=0$ and further there exist terms in the left hand lie in $F^{2}(H)\otimes F^{2}(H)$ while all terms in the right lie in $\delta\otimes\delta$ if $a\neq0$, which is a contradiction. Therefore, there must have $a=0$ and then $k$ is arbitrary.

If $dim\delta=1$, then $H\cong\mathbf{k}[s]$ and $\alpha=s\otimes1-1\otimes s$. From (4.3), one can obtain
\begin{equation*}
  g\otimes g=[g\otimes1+1\otimes g-\Delta(g)]\Delta(g).
\end{equation*}
Comparing the degree of both sides, there must have $g\in \mathbf{k}$.
$\hfill \blacksquare$.
\\

\textbf{Proposition 4.5.} Let $H=U(\delta)$ and $L=He$ be a Lie pseudoalgebra of rank 1 with $[e\ast e]=\alpha\otimes_{H} e$ for some nonzero $\alpha$, then the Reynolds $H$-operator of weight $\lambda$ on $L$ is $R(e)=-\frac{1}{\lambda}e$ or $R(e)=0$ when $\lambda\neq0$; $R(e)=0$ when $\lambda=0$.

\textbf{Proof.} Suppose $R(e)=fe$ for some $f\in H$, then by (2.8), we have
\begin{equation}
  (f\otimes f)\alpha=(1\otimes f+f\otimes1)\alpha\Delta(f)+\lambda(f\otimes f)\alpha\Delta(f).
\end{equation}

If $dim\delta>1$, then $r\neq0$.
Set $deg(f)=d$. If $\lambda=0$, there exists term in the right hand lies in $F^{1}(H)\otimes F^{2d+1}(H)$ while terms in the left hand lie in $F^{d+1}(H)\otimes F^{d+1}(H)$, then there must have $d=0$ and $f\in\mathbf{k}$. Taking this into (4.4), one can obtain $f=0$.
If $\lambda\neq0$, there exists term in the right hand lies in $F^{2d+1}(H)\otimes F^{d+1}(H)$ and terms in the left hand lie in $F^{d+1}(H)\otimes F^{d+1}(H)$, then we also obtain $d=0$ and $f\in\mathbf{k}$. Similarly, one can obtain $f=0$ or $f=-\frac{1}{\lambda}$.

If $dim\delta=1$, then $H\cong\mathbf{k}[s]$ and $\alpha=s\otimes1-1\otimes s$. From (4.4), we have
\begin{equation*}
  f\otimes f=(1\otimes f+f\otimes1)\Delta(f)+\lambda(f\otimes f)\Delta(f).
\end{equation*}
If $\lambda=0$, comparing the degree of both sides, there must have $f\in \mathbf{k}$ and further one can obtain $f=0$. If $\lambda\neq0$, similarly, one can obtain $f\in \mathbf{k}$ and further $f=0$ or $f=-\frac{1}{\lambda}$.
$\hfill \blacksquare$.
\\

\section*{5. Annihilation algebra and Rota-Baxter types operators}
\def\theequation{5. \arabic{equation}}
	\setcounter{equation} {0} \hskip\parindent

Annihilation algebra is an important tool in the study of pseudoalgebras. The annihilation
algebra of the associative pseudoalgebra $CendH = H\otimes H$ is nothing else but the Drinfeld
double (with the obvious comultiplication) of the Hopf algebra $H$. In this section, we present a relationship between Rota-Baxter type $H$-operators on $H$-pseudoalgebra and corresponding operators on their annihilation algebra. During this section, $\mathbf{k}$ will be a field of characteristic 0.\\

Here we recall the definition of the annihilation algebra of an $H$-pseudoalgebra in \cite{BDK}. Let $Y$ be an $H$-bimodule which is a commutative associative $H$-differential algebra both for the left and right action of $H$, for example, $Y=X:=H^{\ast}$.

For a left $H$-module $L$, let $\mathcal{A}_{Y}L=Y\otimes_{H} L$, which has a left action of $H$ in obvious way:
\begin{equation*}
  h(x\otimes_{H} a)=(h\cdot x)\otimes_{H} a, \qquad h\in H, x\in Y, a\in L.
\end{equation*}
If $L$ is an $H$-pseudoalgebra with a pseudoproduct $a\ast b$, we can define a product on $\mathcal{A}_{Y}L$ by:
\begin{equation}
  (x\otimes_{H} a)(y\otimes_{H} b)=\sum_{i} (x\cdot f_{i})(y\cdot g_{i})\otimes_{H} e_{i},
\end{equation}
if $a\ast b=\sum_{i} f_{i}\otimes g_{i}\otimes_{H} e_{i}$.

Let $Y=X$, equation (5.1) endows $\mathcal{A}_{X}L$ an algebra structure, which is called \textbf{annihilation algebra} of $L$.\\

For a given Rota-Baxter type $H$-operator over a finite-rank pseudoalgebra, the following Theorem presents a way to obtain the corresponding Rota-Baxter type operator over its annihilation algebra.\\

\textbf{Theorem 5.1.} Let $H$ be a Hopf algebra and $L$ a free $H$-pseudoalgebra with a basis $\{e_{1}, ...e_{n}\}$. Let $X=H^{\ast}$ and $\mathcal{A}_{X}L$ be the annihilation algebra of $L$, $\xi$ is a left and right $H$-linear map on $X$.

(1). If $\xi(\xi(x)y)=\xi(x)\xi(y)=\xi(x\xi(y))$ and $P$ is an averaging $H$-operator on $L$, then $\xi\otimes_{H} P: \mathcal{A}_{X}L\rightarrow \mathcal{A}_{X}L$ is an averaging operator on $\mathcal{A}_{X}L$;

(2). If $\xi(\xi(x)y)=\xi(x)\xi(y)=\xi(x\xi(y))=\xi^{2}(xy)$ and $P$ is a Nijenhuis $H$-operator on $L$, then $\xi\otimes_{H} P: \mathcal{A}_{X}L\rightarrow \mathcal{A}_{X}L$ is a Nijenhuis operator on $\mathcal{A}_{X}L$;

(3). If $\xi(\xi(x)y)=\xi(x)\xi(y)=\xi(x\xi(y))=\xi(\xi(x)\xi(y))$ and $P$ is a Reynolds $H$-operator on $L$, then $\xi\otimes_{H} P: \mathcal{A}_{X}L\rightarrow \mathcal{A}_{X}L$ is a Reynolds operator on $\mathcal{A}_{X}L$.

\textbf{Proof.} Firstly, we write $(x, y)\leftarrow \alpha =\sum_{i}(x\cdot f_{i})(y\cdot h_{i})\in X$ for convenience if $\alpha=\sum_{i}f_{i}\otimes h_{i}\in H\otimes H$. Since $X$ is right $H$-differential algebra, for any $h\in H$, $\alpha_{1}, \alpha_{2}\in H\otimes H$, by a direct check, we have
\begin{align}
  &(x, y)\leftarrow [\alpha\Delta(h)]=[(x, y)\leftarrow \alpha]\cdot h,\\
  &(x, y)\leftarrow (\alpha_{1}+\alpha_{2})=[(x, y)\leftarrow \alpha_{1}]+[(x, y)\leftarrow \alpha_{2}].
\end{align}

Next, for any $e_{i}, e_{j}\in \{e_{1}, ...e_{n}\}$, we set $e_{i}\ast e_{j}=\alpha_{1}^{ij}\otimes_{H} e_{1}+...+\alpha_{n}^{ij}\otimes_{H} e_{n}$, $P(e_{i})=g_{i1}e_{1}+...+g_{in}e_{n}$ for some $\alpha_{k}^{ij}=\sum_{m}a_{km}^{ij}\otimes b_{km}^{ij}\in H\otimes H$, $g_{it}\in H$ for $i, j, k, t=1, ..., n$.
Then we have
\begin{align*}
  &P(e_{i})\ast P(e_{j})=(g_{i1}\otimes g_{j1}\otimes_{H}1)(e_{1}\ast e_{1})+...+(g_{i1}\otimes g_{jn}\otimes_{H}1)(e_{1}\ast e_{n})+...\\
&+(g_{in}\otimes g_{j1}\otimes_{H}1)(e_{n}\ast e_{1})+...+(g_{in}\otimes g_{jn}\otimes_{H}1)(e_{n}\ast e_{n}) \\
=&\{(g_{i1}\otimes g_{j1})\alpha_{1}^{11}+...+(g_{i1}\otimes g_{jn})\alpha_{1}^{1n}+...+(g_{in}\otimes g_{j1})\alpha_{1}^{n1}+...+(g_{in}\otimes g_{jn})\alpha_{1}^{nn}\}\otimes_{H} e_{1}+...\\
&+\{(g_{i1}\otimes g_{j1})\alpha_{n}^{11}+...+(g_{i1}\otimes g_{jn})\alpha_{n}^{1n}+...+(g_{in}\otimes g_{j1})\alpha_{n}^{n1}+...+(g_{in}\otimes g_{jn})\alpha_{n}^{nn}\}\otimes_{H} e_{n}\\
\overset{\triangle}{=}&\beta_{1}^{ij}\otimes_{H} e_{1}+...+\beta_{n}^{ij}\otimes_{H} e_{n},
\end{align*}
where $\beta_{k}^{ij}=(g_{i1}\otimes g_{j1})\alpha_{k}^{11}+...+(g_{i1}\otimes g_{jn})\alpha_{k}^{1n}+...+(g_{in}\otimes g_{j1})\alpha_{k}^{n1}+...+(g_{in}\otimes g_{jn})\alpha_{k}^{nn}$.
\\
And
\begin{align*}
 &P(e_{i})\ast e_{j}=(g_{i1}\otimes 1\otimes_{H}1)(e_{1}\ast e_{j})+...+(g_{in}\otimes 1\otimes_{H}1)(e_{n}\ast e_{j})\\
=&[(g_{i1}\otimes 1)\alpha_{1}^{1j}+...+(g_{in}\otimes 1)\alpha_{1}^{nj}]\otimes_{H} e_{1}+...+[(g_{i1}\otimes 1)\alpha_{n}^{1j}+...+(g_{in}\otimes 1)\alpha_{n}^{nj}]\otimes_{H} e_{n}\\
=&S_{1}^{ij}\otimes_{H} e_{1}+...+S_{n}^{ij}\otimes_{H} e_{n},
\end{align*}
where $S_{k}^{ij}=(g_{i1}\otimes 1)\alpha_{k}^{1j}+...+(g_{in}\otimes 1)\alpha_{k}^{nj}$.
\begin{align*}
  &e_{i}\ast P(e_{j})=(1\otimes g_{j1}\otimes_{H}1)(e_{i}\ast e_{1})+...+(1\otimes g_{jn}\otimes_{H}1)(e_{i}\ast e_{n})\\
  =&[(1\otimes g_{j1})\alpha_{1}^{i1}+...+(1\otimes g_{jn})\alpha_{1}^{in}]\otimes_{H} e_{1}+...+[(1\otimes g_{j1})\alpha_{n}^{i1}+...+(1\otimes g_{jn})\alpha_{n}^{in}]\otimes_{H} e_{n}\\
  =&T_{1}^{ij}\otimes_{H} e_{1}+...+T_{n}^{ij}\otimes_{H} e_{n},
\end{align*}
where $T_{k}^{ij}=(1\otimes g_{j1})\alpha_{k}^{i1}+...+(1\otimes g_{jn})\alpha_{k}^{in}$.\\

(1). If $P$ is an averaging $H$-operator, by (2.1), we have
\begin{align*}
  &\beta_{1}^{ij}\otimes_{H} e_{1}+...+\beta_{n}^{ij}\otimes_{H} e_{n}\\
  =&(\mathrm{id}\otimes_{H} P)[S_{1}^{ij}\otimes_{H} e_{1}+...+S_{n}^{ij}\otimes_{H} e_{n}]\\
  =&S_{1}^{ij}\otimes_{H} P(e_{1})+...+S_{n}^{ij}\otimes_{H} P(e_{n})\\
  =&S_{1}^{ij}\otimes_{H} (g_{11}e_{1}+...+g_{1n}e_{n})+...+S_{n}^{ij}\otimes_{H} (g_{n1}e_{1}+...+g_{nn}e_{n})\\
  =&[S_{1}^{ij}\Delta(g_{11})+...+S_{n}^{ij}\Delta(g_{n1})]\otimes_{H} e_{1}+...+[S_{1}^{ij}\Delta(g_{1n})+...+S_{n}^{ij}\Delta(g_{nn})]\otimes_{H} e_{n},
\end{align*}
and
\begin{align*}
  &\beta_{1}^{ij}\otimes_{H} e_{1}+...+\beta_{n}^{ij}\otimes_{H} e_{n}\\
  =&(\mathrm{id}\otimes_{H} P)[T_{1}^{ij}\otimes_{H} e_{1}+...+T_{n}^{ij}\otimes_{H} e_{n}]\\
  =&T_{1}^{ij}\otimes_{H} P(e_{1})+...+T_{n}^{ij}\otimes_{H} P(e_{n})\\
  =&T_{1}^{ij}\otimes_{H} (g_{11}e_{1}+...+g_{1n}e_{n})+...+T_{n}^{ij}\otimes_{H} (g_{n1}e_{1}+...+g_{nn}e_{n})\\
  =&[T_{1}^{ij}\Delta(g_{11})+...+T_{n}^{ij}\Delta(g_{n1})]\otimes_{H} e_{1}+...+[T_{1}^{ij}\Delta(g_{1n})+...+T_{n}^{ij}\Delta(g_{nn})]\otimes_{H} e_{n},
\end{align*}
which means
\begin{equation}
  \beta_{k}^{ij}=S_{1}^{ij}\Delta(g_{1k})+...+S_{n}^{ij}\Delta(g_{nk})=T_{1}^{ij}\Delta(g_{1k})+...+T_{n}^{ij}\Delta(g_{nk})\quad for~k=1,...,n.
\end{equation}
Since $\xi$ is right $H$-linear and $\xi(\xi(x)y)=\xi(x)\xi(y)=\xi(x\xi(y))$, we have
\begin{align*}
  (\xi(x), \xi(y))\leftarrow\alpha=&\sum_{i}(\xi(x)\cdot f_{i})(\xi(y)\cdot h_{i})=\sum_{i}\xi(x\cdot f_{i})\xi(y\cdot h_{i})\\=&\xi[\sum_{i}(\xi(x)\cdot f_{i})(y\cdot h_{i})]=\xi[(\xi(x), y)\leftarrow\alpha].
\end{align*}
Similarly, there also have $(\xi(x), \xi(y))\leftarrow\alpha=\xi[(x, \xi(y))\leftarrow\alpha]$.

Therefore, for any $x, y\in X$,
\begin{align*}
  &(\xi(x)\otimes_{H} P(e_{i}))(\xi(y)\otimes_{H} P(e_{j}))=\sum_{k}[(\xi(x), \xi(y))\leftarrow \beta_{k}^{ij}]\otimes_{H}e_{k}\\
  \overset{(5.4)}{=}&\sum_{k}[(\xi(x), \xi(y))\leftarrow [S_{1}^{ij}\Delta(g_{1k})+...+S_{n}^{ij}\Delta(g_{nk})]]\otimes_{H}e_{k}\\
  \overset{(5.3)}{=}&\{(\xi(x), \xi(y))\leftarrow S_{1}^{ij}\Delta(g_{11})+...+(\xi(x), \xi(y))\leftarrow S_{n}^{ij}\Delta(g_{n1})\}\otimes_{H} e_{1}+...\\
  &+\{(\xi(x), \xi(y))\leftarrow S_{1}^{ij}\Delta(g_{1n})+...+(\xi(x), \xi(y))\leftarrow S_{n}^{ij}\Delta(g_{nn})\}\otimes_{H} e_{n}\\
  \overset{(5.2)}{=}&(\xi(x), \xi(y))\leftarrow S_{1}^{ij}\otimes_{H} (g_{11}e_{1}+...+g_{1n}e_{n})+...+(\xi(x), \xi(y))\leftarrow S_{n}^{ij}\otimes_{H} (g_{n1}e_{1}+...+g_{nn}e_{n})\\
  =&[(\xi(x), \xi(y))\leftarrow S_{1}^{ij}\otimes_{H} P(e_{1})]+...
  +[(\xi(x), \xi(y))\leftarrow S_{n}^{ij}\otimes_{H} P(e_{n})]\\
  =&\xi[(\xi(x), y)\leftarrow S_{1}^{ij}]\otimes_{H} P(e_{1})+...+\xi[(\xi(x), y)\leftarrow S_{n}^{ij}]\otimes_{H} P(e_{n})\\
  =&(\xi\otimes_{H} P)\{\underline{[(\xi(x), y)\leftarrow S_{1}^{ij}]\otimes_{H}e_{1}+...+[(\xi(x), y)\leftarrow S_{n}^{ij}]\otimes_{H}e_{n}}\}\\
  =&(\xi\otimes_{H} P)[(\xi(x)\otimes_{H} P(e_{i}))(y\otimes_{H} e_{j})].
\end{align*}
For any $a, b\in L$, if $a=\sum_{i}a_{i}e_{i}$, $b=\sum_{j}b_{j}e_{j}$, then
\begin{align*}
  &(\xi(x)\otimes_{H} P(a))(\xi(y)\otimes_{H} P(b))\\=
  &\sum_{i, j}(\xi(x)\otimes_{H} a_{i}P(e_{i}))(\xi(y)\otimes_{H} b_{j}P(e_{j}))\\
  =&\sum_{i, j}(\xi(x\cdot a_{i})\otimes_{H} P(e_{i})) (\xi(y\cdot b_{j})\otimes_{H} P(e_{j}))\\
  =&\sum_{i, j}(\xi\otimes_{H} P)[(\xi(x\cdot a_{i})\otimes_{H} P(e_{i})) ((y\cdot b_{j})\otimes_{H} e_{j})]\\
  =&(\xi\otimes_{H} P)\{\sum_{i, j}(\xi(x)\otimes_{H} a_{i}P(e_{i}))(y\otimes_{H} b_{j}e_{j})\}\\
  =&(\xi\otimes_{H} P)[(\xi(x)\otimes_{H} P(a))(y\otimes_{H} b)].
\end{align*}
One can also verify
 \begin{equation*}
 (\xi(x)\otimes_{H} P(e_{i}))(\xi(y)\otimes_{H} P(e_{j}))=(\xi\otimes P)[(x\otimes_{H} e_{i})(\xi(y)\otimes_{H} P(e_{j}))]
\end{equation*}
and further
\begin{equation*}
(\xi(x)\otimes_{H} P(a))(\xi(y)\otimes_{H} P(b))=(\xi\otimes_{H} P)[(x\otimes_{H} a)(\xi(y)\otimes_{H} P(b))] \quad for~any~a, b\in L
\end{equation*}
holds in a similar way.
 Thus, $\xi\otimes_{H} P$ is an averaging operator over $\mathcal{A}_{X}L$.
\\

(2). If $P$ is a Nijenhuis $H$-operator, by (2.4), we have
\begin{equation}
  \beta_{k}^{ij}=(S_{1}^{ij}+T_{1}^{ij}-U_{1}^{ij})\Delta(g_{1k})+...+(S_{n}^{ij}+T_{n}^{ij}-U_{n}^{ij})\Delta(g_{nk})\quad for~k=1,...,n,
\end{equation}
where $U_{m}^{ij}=\alpha_{1}^{ij}\Delta(g_{1m})+...+\alpha_{n}^{ij}\Delta(g_{nm})$ for $m=1, 2, ..., n$.

Since $\xi$ is right $H$-linear and $\xi(\xi(x)y)=\xi(x)\xi(y)=\xi(x\xi(y))=\xi^{2}(xy)$, we have
\begin{equation}
(\xi(x), \xi(y))\leftarrow\alpha=\xi[(x, \xi(y))\leftarrow\alpha]=\xi[(\xi(x), y)\leftarrow\alpha]=\xi^{2}[(x, y)\leftarrow\alpha].
\end{equation}
Therefore, for any $x, y\in X$,
\begin{align*}
  &(\xi(x)\otimes_{H} P(e_{i}))(\xi(y)\otimes_{H} P(e_{j}))=\sum_{k}[(\xi(x), \xi(y))\leftarrow \beta_{k}^{ij}]\otimes_{H}e_{k}\\
  \overset{(5.5)}{=}&\sum_{k}[(\xi(x), \xi(y))\leftarrow ((S_{1}^{ij}+T_{1}^{ij}-U_{1}^{ij})\Delta(g_{1k})+...+(S_{n}^{ij}+T_{n}^{ij}-U_{n}^{ij})\Delta(g_{nk}))]\otimes_{H}e_{k}\\
  =&\sum_{k}(\xi(x), \xi(y))\leftarrow (S_{1}^{ij}+T_{1}^{ij}-U_{1}^{ij})\otimes_{H}g_{1k}e_{k}+...+\sum_{k}(\xi(x), \xi(y))\leftarrow (S_{n}^{ij}+T_{n}^{ij}-U_{n}^{ij})\otimes_{H}g_{nk}e_{k}\\
  =&(\xi(x), \xi(y))\leftarrow (S_{1}^{ij}+T_{1}^{ij}-U_{1}^{ij})\otimes_{H}P(e_{1})+...+(\xi(x), \xi(y))\leftarrow (S_{n}^{ij}+T_{n}^{ij}-U_{n}^{ij})\otimes_{H}P(e_{n})\\
  \overset{(5.6)}{=}&\xi[(\xi(x), y)\leftarrow S_{1}^{ij}+(x, \xi(y))\leftarrow T_{1}^{ij}-\xi[(x, y)\leftarrow U_{1}^{ij}]]\otimes_{H}P(e_{1})+...\\
  &+\xi[(\xi(x), y)\leftarrow S_{n}^{ij}+(x, \xi(y))\leftarrow T_{n}^{ij}-\xi[(x, y)\leftarrow U_{n}^{ij}]]\otimes_{H}P(e_{n})\\
  =&(\xi\otimes_{H} P)[\sum_{k}(\xi(x), y)\leftarrow S_{k}^{ij}\otimes_{H}e_{k}+\sum_{k}(x, \xi(y))\leftarrow T_{k}^{ij}\otimes_{H}e_{k}-\sum_{k}\xi[(x, y)\leftarrow U_{k}^{ij}]\otimes_{H}e_{k}]\\
  =&(\xi\otimes_{H} P)[(\xi(x)\otimes_{H} P(e_{i}))(y\otimes_{H} e_{j})+(x\otimes_{H} e_{i})(\xi(y)\otimes_{H} P(e_{j}))-(\xi\otimes_{H} P)[(x\otimes_{H} e_{i})(y\otimes_{H} e_{j})]].
\end{align*}
Similar to (1), one can check
\begin{align*}
&(\xi(x)\otimes_{H}P(a))(\xi(y)\otimes_{H} P(b))\\
=&(\xi\otimes_{H} P)[(\xi(x)\otimes_{H} P(a))(y\otimes_{H} b)+(x\otimes_{H} a)(\xi(y)\otimes_{H} P(b))-(\xi\otimes_{H} P)[(x\otimes_{H} a)(y\otimes_{H} b)]]
\end{align*}
 holds for any $a, b\in L$ and therefore $\xi\otimes P$ is a Nijenhuis operator over $\mathcal{A}_{X}L$.\\

(3). If $P$ is a Reynolds $H$-operator, by (2.8), we obtain
\begin{equation}
  \beta_{k}^{ij}=(S_{1}^{ij}+T_{1}^{ij}+\lambda\beta_{1}^{ij})\Delta(g_{1k})+...+(S_{n}^{ij}+T_{n}^{ij}+\lambda\beta_{n}^{ij})\Delta(g_{nk})\quad for~k=1,...,n.
\end{equation}
Since $\xi$ is right $H$-linear and $\xi(\xi(x)y)=\xi(x)\xi(y)=\xi(x\xi(y))=\xi(\xi(x)\xi(y))$, we have
\begin{equation}
(\xi(x), \xi(y))\leftarrow\alpha=\xi[(x, \xi(y))\leftarrow\alpha]=\xi[(\xi(x), y)\leftarrow\alpha]=\xi[(\xi(x), \xi(y))\leftarrow\alpha].
\end{equation}
Therefore, for any $x, y\in X$,
\begin{align*}
  &(\xi(x)\otimes_{H} P(e_{i}))(\xi(y)\otimes_{H} P(e_{j}))=\sum_{k}[(\xi(x), \xi(y))\leftarrow \beta_{k}^{ij}]\otimes_{H}e_{k}\\
  \overset{(5.7)}{=}&\sum_{k}[(\xi(x), \xi(y))\leftarrow ((S_{1}^{ij}+T_{1}^{ij}+\lambda\beta_{1}^{ij})\Delta(g_{1k})+...+(S_{n}^{ij}+T_{n}^{ij}+\lambda\beta_{n}^{ij})\Delta(g_{nk}))]\otimes_{H}e_{k}\\
  =&\sum_{k}(\xi(x), \xi(y))\leftarrow (S_{1}^{ij}+T_{1}^{ij}+\lambda\beta_{1}^{ij})\otimes_{H}g_{1k}e_{k}+...+\sum_{k}(\xi(x), \xi(y))\leftarrow (S_{n}^{ij}+T_{n}^{ij}+\lambda\beta_{n}^{ij})\otimes_{H}g_{nk}e_{k}\\
  =&(\xi(x), \xi(y))\leftarrow (S_{1}^{ij}+T_{1}^{ij}+\lambda\beta_{1}^{ij})\otimes_{H}P(e_{1})+...+(\xi(x), \xi(y))\leftarrow (S_{n}^{ij}+T_{n}^{ij}+\lambda\beta_{n}^{ij})\otimes_{H}P(e_{n})\\
  \overset{(5.8)}{=}&\xi[(\xi(x), y)\leftarrow S_{1}^{ij}+(x, \xi(y))\leftarrow T_{1}^{ij}+\lambda[(\xi(x), \xi(y))\leftarrow \beta_{1}^{ij}]]\otimes_{H}P(e_{1})+...\\
  &+\xi[(\xi(x), y)\leftarrow S_{n}^{ij}+(x, \xi(y))\leftarrow T_{n}^{ij}+\lambda[(\xi(x), \xi(y))\leftarrow\beta_{n}^{ij}]]\otimes_{H}P(e_{n})\\
  =&(\xi\otimes_{H} P)[\sum_{k}(\xi(x), y)\leftarrow S_{k}^{ij}\otimes_{H}e_{k}+\sum_{k}(x, \xi(y))\leftarrow T_{k}^{ij}\otimes_{H}e_{k}+\lambda\sum_{k}[(\xi(x), \xi(y))\leftarrow \beta_{k}^{ij}]\otimes_{H}e_{k}]\\
  =&(\xi\otimes_{H} P)[(\xi(x)\otimes_{H} P(e_{i}))(y\otimes_{H} e_{j})+(x\otimes_{H} e_{i})(\xi(y)\otimes_{H} P(e_{j}))+\lambda(\xi(x)\otimes_{H} P(e_{i}))(\xi(y)\otimes_{H} P(e_{j}))].
\end{align*}
Similar to (1), one can check
 \begin{align*}
 &(\xi(x)\otimes_{H} P(a))(\xi(y)\otimes_{H} P(b))\\
 =&(\xi\otimes_{H} P)[(\xi(x)\otimes_{H} P(a))(y\otimes_{H} b)+(x\otimes_{H} a)(\xi(y)\otimes_{H} P(b))+\lambda(\xi(x)\otimes_{H} P(a))(\xi(y)\otimes_{H} P(b))]
\end{align*}
holds for any $a, b\in L$ and therefore $\xi\otimes_{H}P$ is a Reynolds operator of weight $\lambda$ on $\mathcal{A}_{X}L$.
$\hfill \blacksquare$.
\\

\textbf{Corollary 5.2} Let $H$ be a Hopf algebra and $L$ a free $H$-pseudoalgebra with a basis $\{e_{1}, ...e_{n}\}$. Let $X=H^{\ast}$ and $\mathcal{A}_{X}L$ be the annihilation algebra of $L$. Define $\xi: X\rightarrow X$ by $\xi(x)=h\cdot x$ where $h\in H$ is a left integral.

(1). If $P$ is an averaging $H$-operator on $L$, then $\xi\otimes_{H} P: \mathcal{A}_{X}L\rightarrow \mathcal{A}_{X}L$ is an averaging operator on $\mathcal{A}_{X}L$;

(2). If $h$ is also a group-like element and $P$ is a Nijenhuis $H$-operator on $L$, then $\xi\otimes_{H} P: \mathcal{A}_{X}L\rightarrow \mathcal{A}_{X}L$ is a Nijenhuis operator on $\mathcal{A}_{X}L$;

(3). If $\epsilon(h)=1$ and $P$ is a Reynolds $H$-operator on $L$, then $\xi\otimes_{H} P: \mathcal{A}_{X}L\rightarrow \mathcal{A}_{X}L$ is a Reynolds operator on $\mathcal{A}_{X}L$.

\textbf{Proof.} $\xi$ is a left and right $H$-linear map since $X$ is $H$-bimodule.

(1). For any $g\in H$, $x, y\in X$, we have
\begin{align*}
  \xi[\xi(x)y](g)=&\langle\xi(x)y, gh\rangle=\langle\xi(x), g_{(1)}h_{(1)}\rangle\langle y, g_{(2)}h_{(2)}\rangle \\
                     =&\langle x, g_{(1)}\underline{h_{(1)}h}\rangle \langle y, g_{(2)}h_{(2)}\rangle=\langle x, g_{(1)}h\rangle \langle y, g_{(2)}h\rangle\\
                     =&\langle\xi(x), g_{(1)}\rangle\langle\xi(y), g_{(2)}\rangle\\
                     =&[\xi(x)\xi(y)](g).
\end{align*}
One can similarly verify that $\xi[x\xi(y)](g)=[\xi(x)\xi(y)](g)$. From Theorem 5.1, $\xi\otimes_{H} P$ is an averaging operator on $\mathcal{A}_{X}L$.

(2). From Theorem 5.1, it is only need to prove $\xi(x)\xi(y)=\xi^{2}(xy)$. Actually, $\forall g\in H$,
\begin{align*}
  \xi^{2}(xy)(g)=&\langle h^{2}\cdot(xy), g\rangle=\langle xy, gh^{2}\rangle=\epsilon(gh)\langle xy, h\rangle\\
  =&\epsilon(g)\langle x, h\rangle\langle y,h\rangle=\langle x, g_{1}h\rangle\langle y, g_{2}h\rangle\\
  =&\langle(h\cdot x)(h\cdot y), g\rangle=[\xi(x)\xi(y)](g).
\end{align*}

(3). It is only need to prove $\xi(x)\xi(y)=\xi[\xi(x)\xi(y)]$. Actually, $\forall g\in H$,
\begin{align*}
  \xi[\xi(x)\xi(y)](g)=&\langle h\cdot(h\cdot x)(h\cdot y), g\rangle=\langle h\cdot x, g_{1}h_{1}\rangle\langle h\cdot y, g_{2}h_{2}\rangle\\
  =&\langle x, g_{1}h_{1}h\rangle\langle y, g_{2}h_{2}h\rangle=\epsilon(gh)\langle x, h\rangle\langle y, h\rangle=\langle x, g_{1}h\rangle\langle y, g_{2}h\rangle\\
  =&\langle (h\cdot x)(h\cdot y), g\rangle=[\xi(x)\xi(y)](g).
\end{align*}
This completes the proof.
$\hfill \blacksquare$.
\\

\section*{6. Rota-Baxter types operators on $H$-conformal algebra}
\def\theequation{6. \arabic{equation}}
	\setcounter{equation} {0} \hskip\parindent

In this section, we recall the notion of Lie $H$-conformal algebra, which is induced by $x$-bracket and first introduced in \cite{BDK}.
Further, we give the definitions of conformal averaging (Nijenhuis, Reynolds) $H$-operator and present a relation between averaging (Nijenhuis, Reynolds) $H$-operator and corresponding conformal $H$-operator.
\\

We start by the $x$-bracket. For any Lie $H$-pseudoalgebra $L$ and $x \in X$, there is a map defined by
\begin{equation*}
  L\otimes L\rightarrow L, \quad a\otimes b\rightarrow [a_{x}b]=\sum_{i}\langle S^{\ast}(x), h_{i}\rangle c_{i},
\end{equation*}
if $[a\ast b]=\sum_{i} h_{i}\otimes 1\otimes_{H} c_{i}$. This map is called an \emph{$x$-brackets}
of $L$. If we write $[a\ast b]=\sum_{i} f_{i}\otimes g_{i}\otimes_{H} c_{i}$, then by (1.7), we have $[a_{x}b]=\sum_{i}\langle S^{\ast}(x), f_{i}S(g_{i(1)})\rangle g_{i(2)}c_{i}$. \\

A general definition of Lie $H$-conformal algebra is given as follow.\\

\textbf{Definition 6.1}(\cite{BDK}) A Lie $H$-conformal algebra is a left $H$-module $L$
equipped with $x$-brackets $[a_{x}b]\in L$ for $a, b\in L, x\in X$,  satisfying the following properties:

\emph{Locality:}
\begin{equation}
  codim\{x\in X\mid [a_{x}b]=0\}<\infty \quad for~any~a, b\in L.
\end{equation}

\emph{$H$-sesqui-linearity:}
\begin{align}
  [ha_{x}b]=&[a_{xh}b],\\
  [a_{x}hb]=&h_{(2)}[a_{S(h_{(1)})x}b].
\end{align}

\emph{Skew-commutativity:}
\begin{equation}
  [a_{x}b]=-\sum_{i}\langle x, S(h_{i(1)})\rangle S(h_{i(2)})[b_{x_{i}}a],
\end{equation}
where $\{h_{i}\}, \{x_{i}\}$ is dual bases in $H$ and $X$.

\emph{Jacobi identity:}
\begin{equation}
  [a_{x}[b_{y}c]]-[b_{y}[a_{x}c]]=[[a_{x_{(2)}}b]_{yx_{(1)}}c].
\end{equation}

By discussion above, if $L$ is a Lie $H$-pseudoalgebra, then there exists a natural Lie $H$-conformal algebra derived from $L$ by letting
\begin{equation}
[a_{x}b]=\sum_{i}\langle S^{\ast}(x), f_{i}S(g_{i(1)})\rangle g_{i(2)}c_{i}
\end{equation}
 if $[a\ast b]=\sum_{i} f_{i}\otimes g_{i}\otimes_{H} c_{i}$. We call $(L, [\cdot_{x}\cdot])$ the \emph{induced Lie $H$-conformal algebra}.
\\

Then we give the definition of conformal averaging (Nijenhuis, Reynolds) $H$-operator, which is both different from averaging (Nijenhuis, Reynolds) $H$-operator on $H$-pseudoalgebra and averaging (Nijenhuis, Reynolds) operator on $\mathbf{k}$-algebra.\\

\textbf{Definition 6.2} Let $L$ be a Lie $H$-conformal algebra,

(1). an $H$-linear map $R$ on $L$ is called a \emph{conformal averaging $H$-operator} if $R[R(a)_{x}b]=R[a_{x}R(b)]=[R(a)_{x}R(b)]$ for any $a, b\in L, x\in X$;

(2). an $H$-linear map $R$ on $L$ is called a \emph{conformal Nijenhuis $H$-operator} if $[R(a)_{x}R(b)]=R[[R(a)_{x}b]+[a_{x}R(b)]-R[a_{x}b]]$ for any $a, b\in L, x\in X$;

(3). an $H$-linear map $R$ on $L$ is called a \emph{conformal Reynolds $H$-operator of weight $\lambda$} if $[R(a)_{x}R(b)]=R[[R(a)_{x}b]+[a_{x}R(b)]+\lambda[R(a)_{x}R(b)]]$ for any $a, b\in L, x\in X$.\\

For an Lie $H$-pseudoalgebra $L$, we affirm that if $R$ is an averaging (respectively, Nijenhuis, Reynolds) $H$-operator on $L$, then $R$ is also a conformal averaging (respectively, Nijenhuis, Reynolds) $H$-operator on its induced Lie $H$-conformal algebra. In order to prove this, we first define a $\mathbf{k}$-linear map $\eta_{x}: H\otimes H\rightarrow H$ by $\eta_{x}(f\otimes g)=\langle S^{\ast}(x), fS(g_{(1)})\rangle g_{(2)}$. Then $\eta_{x}$ has the following properties:
\\

\textbf{Lemma 6.3} Let $H$ be a Hopf algebra and $X=H^{\ast}$, then for any $f, g, h\in H$, $x\in X$,

(1). $\eta_{x}((f\otimes g)\Delta(h))=\eta_{x}(f\otimes g)h$;

(2). $\eta_{x}((h\otimes1)(f\otimes g))=\eta_{x_{(2)}}(h\otimes1)\eta_{x_{(1)}}(f\otimes g)$;

(3). $\eta_{x}((1\otimes h)(f\otimes g))=\eta_{x_{(1)}}(1\otimes h)\eta_{x_{(2)}}(f\otimes g)$.

\textbf{Proof.} By the properties of antipode, we have
\begin{align*}
  \eta_{x}((f\otimes g)\Delta(h))=&\eta_{x}(fh_{(1)}\otimes gh_{(2)})\\
=&\langle S^{\ast}(x), fh_{(1)}S(g_{(1)}h_{(2)(1)})\rangle g_{(2)}h_{(2)(2)}\\
=&\langle S^{\ast}(x), f\underline{h_{(1)}S(h_{(2)(1)})}S(g_{(1)})\rangle g_{(2)}h_{(2)(2)}\\
=&\langle S^{\ast}(x), fS(g_{(1)})\rangle g_{(2)}h\\
=&\eta_{x}(f\otimes g)h.
\end{align*}
\begin{align*}
 \eta_{x}((h\otimes1)(f\otimes g))=&\langle S^{\ast}(x), hfS(g_{(1)})\rangle g_{(2)}\\
 \overset{(1.6)}{=}&\langle S^{\ast}(x)_{(1)}, h\rangle\langle S^{\ast}(x)_{(2)}, fS(g_{(1)})\rangle g_{(2)}\\
 =&\langle S^{\ast}(x_{(2)}), h\rangle\langle S^{\ast}(x_{(1)}), fS(g_{(1)})\rangle g_{(2)}\\
 =&\eta_{x_{(2)}}(h\otimes1)\eta_{x_{(1)}}(f\otimes g).
\end{align*}
\begin{align*}
  \eta_{x}((1\otimes h)(f\otimes g))=&\langle S^{\ast}(x), fS(h_{(1)}g_{(1)})\rangle h_{(2)}g_{(2)}\\
 =&\langle S^{\ast}(x), fS(g_{(1)})S(h_{(1)})\rangle h_{(2)}g_{(2)}\\
 \overset{(1.6)}{=}&\langle S^{\ast}(x)_{(1)}, fS(g_{(1)})\rangle\langle S^{\ast}(x)_{(2)}, S(h_{(1)})\rangle h_{(2)}g_{(2)}\\
 =&\langle S^{\ast}(x_{(1)}), S(h_{(1)})\rangle h_{(2)}\langle S^{\ast}(x_{(2)}), fS(g_{(1)})\rangle g_{(2)}\\
 =&\eta_{x_{(1)}}(1\otimes h)\eta_{x_{(2)}}(f\otimes g).
\end{align*}
This completes the proof.
$\hfill \blacksquare$.
\\

\textbf{Theorem 6.4} Let $L$ be a Lie $H$-pseudoalgebra of rank $n$ and $R$ an averaging (respectively, Nijenhuis, Reynolds) $H$-operator on $L$, then $R$ is also a comformal averaging (respectively, Nijenhuis, Reynolds) $H$-operator on the induced Lie $H$-conformal algebra $L$.

\textbf{Proof.} For any $e_{i}, e_{j}\in \{e_{1}, ...e_{n}\}$, similar to Theorem 5.1, let $e_{i}\ast e_{j}=\alpha_{1}^{ij}\otimes_{H} e_{1}+...+\alpha_{n}^{ij}\otimes_{H} e_{n}$ and $R(e_{i})=g_{i1}e_{1}+...+g_{in}e_{n}$ for some $\alpha_{k}^{ij}=\sum_{m}a_{km}^{ij}\otimes b_{km}^{ij}\in H\otimes H$, $g_{it}\in H$ for $i, j, k, t=1, ..., n$. We set
\begin{align*}
  &[R(e_{i})\ast e_{j}]=S_{1}^{ij}\otimes_{H} e_{1}+...+S_{n}^{ij}\otimes_{H} e_{n},\\
  &[e_{i}\ast R(e_{j})]=T_{1}^{ij}\otimes_{H} e_{1}+...+T_{n}^{ij}\otimes_{H} e_{n},\\
  &[R(e_{i})\ast R(e_{j})]=\beta_{1}^{ij}\otimes_{H} e_{1}+...+\beta_{n}^{ij}\otimes_{H} e_{n},
\end{align*}
for some $S_{k}^{ij}, T_{k}^{ij}, \beta_{k}^{ij}\in H\otimes H$.

If $R$ is an averaging $H$-operator on $L$, then by (2.1), we have
\begin{align*}
  \beta_{k}^{ij}=&S_{1}^{ij}\Delta(g_{1k})+...+S_{n}^{ij}\Delta(g_{nk})\\
  =&T_{1}^{ij}\Delta(g_{1k})+...+T_{n}^{ij}\Delta(g_{nk})
\end{align*}
for $k=1, 2, ...n$.

Thus, for any $x\in X$, by Lemma 6.3, we have
\begin{align*}
  [R(e_{i})_{x}R(e_{j})]=&\sum_{k=1}^{n}\eta_{x}(\beta_{k}^{ij})e_{k}=\sum_{k=1}^{n}\eta_{x}(S_{1}^{ij}\Delta(g_{1k})+...+S_{n}^{ij}\Delta(g_{nk}))e_{k}\\
  =&\sum_{k}[\eta_{x}(S_{1}^{ij})g_{1k}e_{k}+...+\eta_{x}(S_{n}^{ij})g_{nk}e_{k}]\\
=&\eta_{x}(S_{1}^{ij})(g_{11}e_{1}+...+g_{1n}e_{n})+...+\eta_{x}(S_{n}^{ij})(g_{n1}e_{1}+...+g_{nn}e_{n})\\
=&\eta_{x}(S_{1}^{ij})R(e_{1})+...+\eta_{x}(S_{n}^{ij})R(e_{n})\\
=&R[\sum_{k}\eta_{x}(S_{k}^{ij})e_{k}]\\
=&R[R(e_{i})_{x}e_{j}].
\end{align*}
For any $a=\sum_{i}h_{i}e_{i}, b=\sum_{i}f_{i}e_{i}\in L$,
\begin{align*}
  [R(a)_{x}R(b)]=&\sum_{i,j}[h_{i}R(e_{i})_{x}f_{j}R(e_{j})]\\
  \overset{(6.2),(6.3)}{=}&\sum_{i,j}f_{j(2)}[R(e_{i})_{S(f_{j(1)})xh_{i}}R(e_{j})]\\
  =&\sum_{i,j}f_{j(2)}R[R(e_{i})_{S(f_{j(1)})xh_{i}}e_{j}]\\
  =&\sum_{i,j}R[h_{i}R(e_{i})_{x}f_{j}e_{j}]\\
  =&R[R(a)_{x}b].
\end{align*}
One can also obtain $[R(a)_{x}R(b)]=R[a_{x}R(b)]$ in a similar way. Thus, $R$ is a comformal averaging $H$-operator on the induced Lie $H$-conformal algebra.\\

If $R$ is a Nijenhuis $H$-operator on $L$, then from (2.4), we have
\begin{equation*}
  \beta_{k}^{ij}=(S_{1}^{ij}+T_{1}^{ij}-U_{1}^{ij})\Delta(g_{1k})+...+(S_{n}^{ij}+T_{n}^{ij}-U_{n}^{ij})\Delta(g_{nk})
\end{equation*}
for $k=1, 2, ...n$, where $U_{m}^{ij}=\alpha_{1}^{ij}\Delta(g_{1m})+...+\alpha_{n}^{ij}\Delta(g_{nm})$. By Lemma 6.3, we obtain
\begin{equation*}
\eta_{x}(U_{m}^{ij})=\sum_{l}\eta_{x}(\alpha_{l}^{ij})g_{lm}
\end{equation*}
for any $x\in X$. Then
\begin{align*}
  [R(e_{i})_{x}R(e_{j})]=&\sum_{k=1}^{n}\eta_{x}(\beta_{k}^{ij})e_{k}=\sum_{k=1}^{n}\eta_{x}[(S_{1}^{ij}+T_{1}^{ij}-U_{1}^{ij})\Delta(g_{1k})+...
  +(S_{n}^{ij}+T_{n}^{ij}-U_{n}^{ij})\Delta(g_{nk})]e_{k}\\
  =&\sum_{k=1}^{n}\eta_{x}(S_{1}^{ij}+T_{1}^{ij}-U_{1}^{ij})g_{1k}e_{k}+...+\sum_{k=1}^{n}\eta_{x}(S_{n}^{ij}+T_{n}^{ij}-U_{n}^{ij})g_{nk}e_{k}\\
  =&\eta_{x}(S_{1}^{ij}+T_{1}^{ij}-U_{1}^{ij})R(e_{1})+...+\eta_{x}(S_{n}^{ij}+T_{n}^{ij}-U_{n}^{ij})R(e_{n})\\
  =&\sum_{k=1}^{n}\eta_{x}(S_{k}^{ij})R(e_{k})+\sum_{k=1}^{n}\eta_{x}(T_{k}^{ij})R(e_{k})-\sum_{k=1}^{n}\eta_{x}(U_{k}^{ij})R(e_{k})\\
  =&R[[R(e_{i})_{x}e_{j}]+[(e_{i})_{x}R(e_{j})]]-R[\sum_{k,l}\eta_{x}(\alpha_{l}^{ij})g_{lk}e_{k}]\\
  =&R[[R(e_{i})_{x}e_{j}]+[(e_{i})_{x}R(e_{j})]]-R(\sum_{l}\eta_{x}(\alpha_{l}^{ij})R(e_{l}))\\
  =&R[[R(e_{i})_{x}e_{j}]+[(e_{i})_{x}R(e_{j})]-R[(e_{i})_{x}e_{j}]].
\end{align*}
Similar to discussion above, one can obtain $[R(a)_{x}R(b)]=R[[R(a)_{x}b]+[a_{x}R(b)]-R[a_{x}b]]$ for any $a, b\in L$. Thus, $R$ is a comformal Nijenhuis $H$-operator on the induced Lie $H$-conformal algebra.\\

If $R$ is a Reynolds $H$-operator, then from (2.8), we have
\begin{equation*}
  \beta_{k}^{ij}=(S_{1}^{ij}+T_{1}^{ij}+\lambda\beta_{1}^{ij})\Delta(g_{1k})+...+(S_{n}^{ij}+T_{n}^{ij}+\lambda\beta_{n}^{ij})\Delta(g_{nk})
\end{equation*}
for $k=1, 2, ...n$.

Then for any $x\in X$,
\begin{align*}
  [R(e_{i})_{x}R(e_{j})]=&\sum_{k=1}^{n}\eta_{x}(\beta_{k}^{ij})e_{k}=\sum_{k=1}^{n}\eta_{x}[(S_{1}^{ij}+T_{1}^{ij}+\lambda\beta_{1}^{ij})\Delta(g_{1k})+...
  +(S_{n}^{ij}+T_{n}^{ij}+\lambda\beta_{n}^{ij})\Delta(g_{nk})]e_{k}\\
  =&\eta_{x}(S_{1}^{ij}+T_{1}^{ij}+\lambda\beta_{1}^{ij})R(e_{1})+...+\eta_{x}(S_{n}^{ij}+T_{n}^{ij}+\lambda\beta_{n}^{ij})R(e_{n})\\
  =&\sum_{k=1}^{n}\eta_{x}(S_{k}^{ij})R(e_{k})+\sum_{k=1}^{n}\eta_{x}(T_{k}^{ij})R(e_{k})+\lambda\sum_{k=1}^{n}\eta_{x}(\beta_{k}^{ij})R(e_{k})\\
  =&R[[R(e_{i})_{x}e_{j}]+[(e_{i})_{x}R(e_{j})]+\lambda[R(e_{i})_{x}R(e_{j})]]
\end{align*}
Similarly, one can obtain $[R(a)_{x}R(b)]=R[[R(a)_{x}b]+[a_{x}R(b)]+\lambda[R(a)_{x}R(b)]]$ for any $a, b\in L$. Thus, $R$ is a comformal Reynolds $H$-operator of weight $\lambda$ on the induced Lie $H$-conformal algebra.
$\hfill \blacksquare$.
\\

\section*{Acknowledgements}

The second author thanks the financial support of the National Natural Science Foundation of
 China (Grant No. 12271089).

\end{document}